\documentclass[a4paper,12pt]{article}

\usepackage{a4wide}
\usepackage{graphicx,subfigure,hhline}
\usepackage[colorlinks=true,linkcolor=blue,citecolor=red]{hyperref}
\usepackage{xcolor}
\usepackage{amssymb,amsthm,empheq,bbold}

\usepackage[ruled,vlined,linesnumbered]{algorithm2e}
\usepackage{caption} 
\usepackage{enumerate}
\usepackage{url}
\usepackage{multirow, filecontents}
\usepackage{array}
\usepackage{xcolor}
\usepackage{calc,mathtools}
\usepackage{float}
\usepackage{babel}

\usepackage{cancel}
\usepackage[normalem]{ulem}

\newtheorem{remark}{Remark}

\newtheorem{lemma}{Lemma}

\newcommand{\be}{\begin{equation}}
\newcommand{\ee}{\end{equation}}
\newcommand{\nn}{\nonumber}

\newcommand{\dis}{\displaystyle}
\newcommand{\pa}{\partial}
\newcommand{\bv}{\mathbf{v}}

\newcommand{\bx}{\mathbf{x}}
\newcommand{\bvp}{\mathbf{v'}}

\newcommand{\fvec}{\underline{\mathbf{f}}}

\begin{document}
	
	\title{Reaction--diffusion systems derived from \\ kinetic models for Multiple Sclerosis}
	
	\author{Romina Travaglini {$^{1*}$}, Jo\~ao Miguel Oliveira{$^{1**}$}\\[1em]
		{\footnotesize $^1$ Centre of Mathematics, University of Minho,}\\
{\footnotesize  Campus of Gualtar, Braga, Portugal, 4710-057}
\\[0.1em]
{* \footnotesize b13201@math.uminho.pt}
		\\[0.5em]
{**\footnotesize b6885@math.uminho.pt}}

	\date{}
	\maketitle

	\begin{abstract}
		We present a mathematical study for the development of Multiple Sclerosis in which a spatio-temporal kinetic model describes, at mesoscopic level, the dynamics  of a high number of interacting agents. We consider both interactions among different populations of human cells and motion of immune cells, stimulated by cytokines. Moreover, we reproduce the consumption of myelin sheath due to anomalously activated lymphocytes and its restoration by oligodendrocytes. Successively, we fix a small time parameter and assume that the considered processes occur at different scales. This allows to perform a formal limit, obtaining macroscopic reaction-diffusion equations for the number densities with a chemotaxis term. A natural step is then to  study the system, inquiring about the formation of spatial patterns through a Turing instability analysis of the problem and basing the discussion on microscopic parameters of the model. In particular, we get spatial patterns oscillating in time that may reproduce brain lesions characteristic of different phases of the pathology.
	\end{abstract}

	\smallskip
	
	
	\smallskip
	
	
\vspace{1cm}
\emph {Declaration of interest:} The authors declare that they have no known competing financial interests or personal relationships that could have appeared to influence the work reported in this paper.

	\section{Introduction}
Multiple Sclerosis is one of the most severe inflammatory diseases of the central nervous system. Early observations of brain damage due to Multiple Sclerosis date back to the late 19th and early 20th centuries. Since then, it has been evident that a key pathological feature is the formation of plaques in the white matter. These plaques represent lesions driven by inflammation of the myelin sheath in the brain, essential for facilitating the transmission of cerebral impulses.

The most accepted mechanism underlying Multiple Sclerosis is an autoimmune response. More precisely, immune system cells, that can be activated by self-antigens, are not normally suppressed by regulatory cells. Instead, these immune cells may erroneously attack healthy tissue, initiating an autoimmune cascade. For the particular case of Multiple Sclerosis considered here, T-cells, B-cells, macrophages, and microglia become activated against antigens expressed in myelin and oligodendrocytes (cells responsible for myelin production).

From the first studies\cite{babinski1885recherches} it became clear that the primary characteristic of Multiple Sclerosis is the inflammatory process leading to myelin injury. However, demyelination patterns have been observed to be, at an early stage of the disease, homogeneous at the individual level but heterogeneous between different patients\cite{lucchinetti2000heterogeneity}. This suggests a wide range of immune mechanisms underlying the formation of plaques. Moreover, the clinical course of Multiple Sclerosis and features of lesions, as well as associated irreversible neurological symptoms, are extremely varied between patients. Different stages of the disease, including the possibility of plaque restoration, contribute to this complexity.

In this framework, we provide a mathematical model simulating the appearance of myelin lesions along with the remyelination process. Several mathematical models have been carried out to study Multiple Sclerosis. A particularly detailed system of partial differential equations for cells and substances involved in plaque formation can be found in  \cite{moise2021mathematical}. Simpler ordinary differential equation systems have been proposed, instead, to describe the relapsing-remitting dynamics\cite{elettreby2020simple,frascoli2022dynamics} and brain damage\cite{kotelnikova2017dynamics}. Models able to reproduce the formation of lesions via Turing instability have also been proposed, for example, in  \cite{lombardo2017demyelination} 
for the particular case of Balo's Sclerosis.

Models cited above, though, are straightforwardly provided for macroscopic densities of cells and substances. All the parameters in the equations come from experimental observation or heuristic considerations. Nevertheless, the microscopic interplay among cells and molecules is at the base of the problem considered. Consequently, a mathematical description of microscopic dynamics would be desirable. Such a description would lead to the derivation of a coherent macroscopic scenario for observable phenomena while maintaining a close relationship with micro-scale features.

{The cellular dynamics of immune response is an example of a complex system, composed of many heterogeneous living individuals. These individuals interact stochastically within themselves and with the external environment. To study such phenomena, the kinetic theory of active particles\cite{survey} proves to be an apt tool. It draws inspiration from the mathematical kinetic theory of inert matter, especially from the Boltzmann equation. At the same time, though, it introduces methods that recognize the disparities between physical inert systems and living organisms. The most relevant distinction is that the Boltzmann equation exploits the assumption of a rarefied flow. In this case, interactions are exclusively binary and occur at short range. In contrast, living individuals interact through mechanisms involving sensitivity and visibility, related to non-locality and multiple interplays. Additionally, encounters among classical particles preserve mechanical quantities, while such conservations are often lost in living systems. 
	
	First examples of kinetic {theory} modeling of living individuals, based on systems of integro-differential equations, can be found in  \cite{bellomo1994dynamics}. Successively, they have been applied to many physical problems (as listed in  \cite{survey})}.
For the autoimmune case, though, the literature is not very extensive. Some results are presented in  \cite{delitala2013mathematical,kolev2018mathematical}, followed by  \cite{della2022mathematical,ramos2019kinetic}. In the last cited works, the authors apply the kinetic theory for active particles to populations of self-antigen presenting cells, self-reactive T cells, and immunosuppressive cells. Each one of these populations is endowed with a microscopic functional state.
In these cases, a proper integration leads to a macroscopic description for the dynamics of biologically relevant quantities over time. However, a description also in terms of spatial variables seems suitable. Specifically, it would reproduce the movement of cells along biological tissues. Moreover, it is widely accepted that leukocyte migration (related to inflammation) is regulated by chemotactic motion induced by cytokines\cite{taub1995chemotaxis}.

For this reason, we propose a kinetic { theory} description allowing the derivation of partial-differential equations of reaction-diffusion type. The derivation is performed through a time scaling, followed by a proper diffusive limit. {This tool has been applied in different areas, from classical Boltzmann theory of gas dynamics\cite{bisi2006reactive,bisi2022reaction,lachowicz2002microscopic} to the description of cells and tissues (see  \cite{burini2019multiscale} and references therein).}
The chemotactic term will be derived following the procedure outlined in  \cite{alt1980biased}. In this work, a bias coming from an external field, and influencing the run-and-tumble movement of cells is introduced through a turning operator.  {In addition, the bias will be interpreted as a first-order perturbation of a symmetric probability of the velocity, as done in  \cite{othmer1988models}. More precisely, the perturbation results from the gradient of a chemotactic attracting substance.  It is worth pointing out here that such derivation has been {thoroughly} refined in  \cite{othmer2000diffusion,othmer2002diffusion}, {while also being} applied to the general study of cancer\cite{bellouquid2006mathematical}}. { More recently, further modeling features have been included in the run and tumbling dynamics such as, for example, non-local bias\cite{loy2020kinetic}, cell adhesion\cite{loy2020modelling} or motion in a confined environment\cite{fu2023confined}.} Finally, {such an approach has been applied to autoimmune diseases. In  \cite{immuneProc}, indeed, the model proposed in  \cite{della2022mathematical} has been extended including the chemotactic motion of T-cells driven by cytokines.}

The work is structured as follows. {In Section \ref{Sec2}, a {detailed} biological description of the model is provided. In Section \ref{Sec3}, we build up the kinetic {{theory}} structure of the dynamics. Then, in Section \ref{Sec4},} we propose a time scaling of the model leading to a macroscopic system of reaction-diffusion equations with a chemotaxis term. The system obtained is analyzed in Section \ref{Sec5} from the Turing instability point of view, finding conditions on parameters that induce the formation of patterns. Such patterns reproduce the demyelinating process characteristic of Multiple Sclerosis, and the depiction of different phases listed above is obtained by means of numerical simulations. Finally, in Section \ref{Sec6} some concluding remarks and future perspectives are mentioned.

{\section{Biological setting} \label{Sec2}
	This section serves as a preliminary overview, providing the essential biological context to justify our model and describe the processes involved in Multiple Sclerosis.
	
	Immunity, and autoimmunity, in particular, is a biological mechanism involving a huge number of cells interacting through { many} complicated processes of activation and migration. These processes are predominantly regulated by specific molecules such as antigens and cytokines\cite{akkaya2019regulatory}.
	{The underlying dynamics in autoimmunity originate when CD4$^+$- T cells present in the periphery are stimulated to differentiate into two particular subgroups (Th1 and Th17).  These groups produce {a} type of proinflammatory cytokines (chemokines), i.e. molecules performing as attractants for other cells. Chemokines may enter the central nervous system through the blood-brain barrier and recruit other immune cells like B cells, T cells, and macrophages. These cells are subsequently activated by cognate antigen\cite{dhaiban2021role} carried by self-antigen presenting cells as macrophages, dendritic cells, and B cells\cite{pone2010toll}. The importance of chemokines and cytokines inducing self-reactive cell migration, through chemotactic motion, has been confirmed by several studies\cite{luster1998chemokines}. {Additionally,} proinflammatory cytokines stimulate clonal expansion of CD8$^+$ T-cell and B-cells, promoting the immune cascade.}
	
	It is necessary to mention that self-reactive immune cells may be present in peripheral tissues also in non-pathological conditions\cite{danke2004autoreactive}. In this case, though,  specific cell populations act as immunosuppressors. These can be represented by regulatory T lymphocytes (Tregs)  and natural killer cells that may inhibit or eliminate self-antigen presenting cells and self-reacting cells. {In Multiple Sclerosis patients, the efficiency of natural killers and Tregs to cause lysis of dendrocytes and to kill activated T-cells is lacking} (for a deeper overview of immunosuppressive cells in Multiple Sclerosis, we address the reader to  \cite{hoglund2013one,mimpen2020natural,zozulya2008role} and references therein).
	
	We point out here that the primary focus of this work is to prioritize the overall process of myelin destruction and restoration. Therefore, we will avoid an overly nuanced differentiation in all populations involved, recognizing the potential for more detailed models in future work. Considering the model proposed in  \cite{della2022mathematical} and extended in  \cite{immuneProc}, {to adapt} it to the Multiple Sclerosis {condition}, we are going to {focus} our analysis upon the following populations of cells and biological substances: self-antigen presenting cells, self-reactive leukocytes, immunosuppressive cells, cytokines, and myelin. 
	
	{Let us consider the effects of the autoimmune cascade on the central nervous system. Active lesions, characterized by acute inflammation and blood-brain barrier damage, are more prevalent during the early phase of the disease. Simultaneously, remyelination processes occur at this stage. This is due to the action of oligodendrocytes and results in the formation of ``shadow" plaques. 
		
		Evidences of this demyelination-remyelination phenomenon were extensively studied by John Prineas in his works from the Eighties and Nineties.  Most frequently,  patients first go through a relapsing-remitting stage (RRMS), in which extensive remyelination is more evident. This stage may last for several years, but it is usually followed by a secondary progressive phase (SPMS). During this second phase, active myelin lesions and blood-brain barrier injuries are less frequent. At this stage, there is evidence of other processes causing demyelination and neurodegeneration. In certain cases, instead, the relapsing-remitting phase may not occur. In such circumstances, an immediate unremitting process leads to the so-called primary progressive disease (PPMS). {For a more complete overview of medical studies on Multiple Sclerosis we address the reader to  \cite{lassmann2005multiple,lassmann2007immunopathology,lassmann2012progressive,mahad2015pathological} and references therein.}}
	Also, the different portions of areas undergoing demyelination-remyelination, in diverse phases of the disease, have been extensively studied, see\cite{lassmann2007immunopathology} and reference therein.
	
	Characterizing the timing of lesion appearance and their subsequent restoration remains a difficult task. Some studies\cite{filippi1998magnetization} suggest that focal changes in the white matter can be observed weeks before the actual development of lesions in the same site. Further research \cite{zhao1996mri} shows that, in general, new active lesions enlarge until reaching a maximum size in  2-4 weeks. Afterward, they can stay stable for long periods or be restored during the following 4-8 weeks. For this reason, we will consider two distinct phases for both the consumption and restoration of myelin. We will also distinguish the different evolution times for each phase.
	
	{Now, we outline  processes among those described above that will be  integrated into the mathematical model:
		{	\begin{itemize}
				\item Interplay of activation-proliferation-regulation. Self-reactive leukocytes are activated through interaction with self-antigen presenting cells. In addition, when interacting, self-reactive leukocytes, self-antigen presenting cells, and immunosuppressive cells may be stimulated to proliferate. Finally, immunosuppressive cells can induce apoptosis (programmed cell death) of self-reactive leukocytes and self-antigen presenting cells.
				\item Production of cytokines. When activated by self-antigen presenting cells, self-reactive leukocytes produce particular molecules called cytokines. Cytokines play a fundamental role in the autoimmune cascade.
				\item Migration. Self-reactive leukocytes are driven by the chemical gradient of cytokines. In other words, they tend to migrate towards sites where cytokines concentration is higher.
				\item Consumption and restoration. Myelin sheath is attacked and consumed by activated self-reactive leukocytes. Afterward, its restoration is performed by oligodendrocytes (cells responsible for myelin production).
		\end{itemize}}
		A thorough analysis and discussion of the involved parameters will be included in the mathematical description of these dynamics. Our goal is to develop a model that can reproduce macroscopic outcomes comparable to the diverse scenarios associated with the phases of Multiple Sclerosis.

		

		\section{{ Mathematical description via kinetic theory of active particles}}
		\label{Sec3}

		In this section, we propose a mathematical description of biological dynamics that are relevant to our aim of reproducing the formation and restoration of myelin lesions. 
		
		{Processes described in the previous section may be well {represented} by means of {the} kinetic theory of active particles\cite{survey}.} For the dynamics involving cellular populations and cytokines, we shall refer to the model proposed in  \cite{immuneProc}, considering the larger set of self-reactive leukocytes instead of only self-reactive T-cells.
		{We} will also propose a more general form for the diffusion and chemotaxis operator for self-reactive leukocytes. 
		Then, we will introduce {here} the description of the demyelination-remyelination process. In particular, we shall distinguish three states for the myelin (sane, partially sane, destroyed). At the first stage, we write {the}  equations for the {dynamics of the} first two states. 
		
		The behavior of each {population} involved {in the description} will be {specified} by its distribution function. 
		More precisely, all distributions will depend on time $t\in \mathbb R_0^+$ and space
		$\bx \in \Gamma_{\bx}$, with $\Gamma_{\bx}$ a bounded {domain} of $\mathbb R^n$ (we are going to perform our analysis in general for any space dimension $n\,=\,1,\,2,\,3$). {Additionally, being cells regarded as active particles, we introduce the activity variable $u\in[0,1]$. In the framework of the kinetic theory of active particles, cell populations represent functional subsystems. Each subsystem includes all the interacting agents sharing the same task, and the activity variable is the level of activation required to perform that task.} Distribution functions of self-reactive leukocytes and cytokines will also depend on the velocity variable.
		{In fact, we are going to take into account the spatial diffusion of these quantities, along with the chemotaxis interplay.} More precisely, we self-reactive leukocytes velocity $\bv\in \Gamma_{R}\,=\,\mathcal V\mathbb B^n$ and cytokines velocity  $\bv\in \Gamma_{C}\,=\,\mathcal W \mathbb B^n$, being $\mathcal V$ and $\mathcal W$ the maximal speed of self-reactive leukocytes and cytokines, respectively, and $\mathbb {B}^{n}$ the unit ball in $\mathbb{R}^n$.
		{ Functional subsystems of the model, along with their distributions function{s,} are then listed below.}
		
		\begin{itemize}
			\item[] {\large $A$} -- Self-antigen presenting cells,
			distribution function $\,f_A(t,\bx,u)$, activity: ability to activate self-reactive leukocytes.
			\vskip2mm
			\item[] {\large $S$} -- Immunosuppressive cells,
			distribution function $\,f_S(t,\bx,u)$, activity:  ability to reduce the activity or suppress self-antigen presenting cells and self-reactive leukocytes.
			\vskip2mm
			\item[] {\large $R$} -- Self-reactive leukocytes,
			distribution function $\,f_R(t,\bx,\bv,u)$, activity: production of cytokines.
			\vskip2mm
			\item[] {\large $C$} -- Cytokines, distribution function $\,f_C(t,\bx,\bv)$.
			\vskip2mm
			\item[] {\large $E_1$} -- Sane myelin, distribution function $\,f_{E_1}(t,\bx)$.
			\vskip2mm
			\item[] {\large $E_2$} -- Partially sane myelin, distribution function $\,f_{E_2}(t,\bx)$.
			\vskip2mm
			\item[] {\large $E $} -- Destroyed myelin, distribution function $\,f_{E}(t,\bx)$.
			\vskip2mm
			{\item[] {\large $O $} -- Oligodendrocytes, that we consider a constant population.}
		\end{itemize}
		In this work, besides the cellular dynamics, we are also interested in the time-space evolution of the populations listed above. Therefore, we define the corresponding macroscopic densities. For cell populations $A$, $R$, $S$, and cytokines $C$, we shall perform  integral of the distribution functions over the space of remaining variables (activity, velocity, or both)
		\be \nn
		\begin{array}{c}
			\dis A(t,\bx) \,=\, \int_0^1 \,f_A(t, \bx, u)\,du,
			\\
			\\
			\dis R(t,\bx)\,=\,\int_0^1 \rho_R(t, \bx, u)\,du, \quad
			\mbox{with} \quad \rho_R(t, \bx, u)\,=\,\int_{\Gamma_R} \,f_R(t,\bx,\bv,u)\,d\bv,
			\\
			\\
			\dis S(t,\bx)\,=\,\int_0^1 \,f_s(t, \bx, u)\,du, 
			\\
			\\
			\dis C(t, \bx)\,=\,\int_{\Gamma_C} \,f_C(t,\bx,\bv)\,d\bv.
		\end{array}
		\ee
		Instead, for the different stages of myelin, we are simply going to change the notation from  
		$\,f_{E_1}(t,\bx)$, $\,f_{E_2}(t,\bx)$, $\,f_{E}(t,\bx)$ to ${E_1}(t,\bx)$,  ${E_2}(t,\bx)$,  ${E}(t,\bx)$, respectively, when dealing with macroscopic quantities.

		The evolution of each distribution functio{n} for $A$, $S$, $E_1$, $E_2$ and $E$ is given by an integro-differential equation of the form
		\be
		\dis\frac{\pa \,f_I}{\pa t}\,=\,  \mathcal N_I(\fvec)+ \mathcal I_I(\fvec),\quad I\,=\,A,\,S,\,E_1,\,E_2,\,E.
		\label{eq:kin1}
		\ee
		The evolution of distribution functions of $R$ and $C$ is given by
		{a similar equation containing also a drift term }
		\be
		\dis \frac{\pa \,f_I}{\pa t}+\bv\cdot\nabla_{\bx}\, \,f_I\,=\, \mathcal {L}_I(f_I)+ \mathcal N_I(\fvec)+ \mathcal I_I(\fvec), \quad I\,=\,R,\,C,
		\label{eq:kin2}
		\ee
		being $\fvec$ the vector of all distribution functions. 
		The terms of type $\mathcal N_I$ include integral operators describing the outcomes of interactions among agents, which may be conservative, proliferative, or destructive. Besides, terms of type $\mathcal I_I$ incorporate operators related to proliferation or destruction based on other processes. Finally, terms of type $\mathcal L_I$ involve turning operators concerning the movement of cells or molecules in the spatial domain. The complete set of processes, along with their corresponding operators, will be outlined in the following.
		
		\subsection{{ Interaction} operators}
		
		The biological dynamics involving activation, proliferation, suppression of cells, and cytokines production, can be seen as a result of { interactions} among self-antigen presenting cells, self-reactive leukocytes, and 
		immunosuppressive cells. More precisely:
		\begin{itemize}
			\item { Interactions} between $A$ and $R$ can be conservative (with the only result of activity $u$ increasing for both populations) or proliferative (with the birth of a new cell) either for $A$ or for {$R$}. In the proliferative case, the newborn cell has the same activity as its mother cell. In addition, { interactions} between $A$ and $R$ may also result in the production of cytokines by leukocytes.
			\item { Interactions} between $A$ and $S$ can be  conservative (with the activity of $A$ decreasing), destructive for $A$, {i.e.} the $S$ cell causes apoptosis of the $A$ cell, or proliferative for $S$. In this case, the newborn cell has the same activity as its mother cell.
			\item Interactions between $R$ and $S$ are conservative (in which the activity of $R$ decreases) or destructive, that result in the death of self-reactive leukocytes.
		\end{itemize}
		For a more detailed explanation of biological justifications for the choices performed above, we address the reader to  \cite{ramos2019kinetic}.
		As mentioned, self-activated leukocytes $R$ attack the sane myelin, and we {describe the primary damaging interaction as}
		\be \label{dem1} E_1+R\rightarrow E_2+R,\ee and the total destruction, expressed by the interaction
		\be \label{dem2} E_2+R\rightarrow E+R.\ee
		We are now able to write the explicit formulae for { interaction} operators, reporting the ones introduced in  \cite{della2022mathematical,ramos2019kinetic} for populations $A$, $R$, $S$, $C$. {We omit the operators relevant to the conservative processes, since they do not give contribution in the macroscopic description}. We also introduce the new ones for $E_1$ and $E_2$ (writing the operators for processes   \eqref{dem1} and \eqref{dem2} separately). 
		
		{We assume that, in a mean-field-like approximation, the whole density contributes to the outcomes of interaction. Indeed, there is no {distinction} among different values of the activity variable of interacting agents. Due to this, the transition probabilities, 
			{typical} of the kinetic theory of active particles, reduce to constant transition rates.} We shall adopt the notation $c_{IJ}$, $p_{IJ}$, $d_{IJ}$, $b_{IJ}$ for the conservative, proliferative, destructive and myelin-damaging rates, respectively, where indices $I, J$ vary in the set $\{1,\ldots,6\}$ corresponding to quantities $A,\, R,\, S,\, C,\,E_1,\,E_2$.
		Thus we have 
		\be
		\begin{aligned} 
			\mathcal N_A(\fvec)&\,=\,
			- c_{12} (u -1)^2 \,f_A(t,\bx,u)  \int_{\Gamma_{R}}   \int_{0 }^1  \,f_R (t,\bx,\bv,w) d w \,d\bv \\ 
			&\qquad +p_{12}  \,f_A(t,\bx,u) \int_{0 }^1 \,f_R(t,\bx,\bv,w) d w   - d_{13} \,f_A(t,\bx,u)  \int_{0 }^1  \,f_S (t,\bx,w) d w,
			\\[2mm]
			\mathcal N_S(\fvec)&\,=\,p_{13}  \,f_S(t,\bx,u) \int_{0 }^1 \,f_A(t,\bx,w) d w,
			\\[2mm]
			\mathcal N_R(\fvec)&\,=\,
			p_{21}  \,f_R(t,\bx,\bv,u) \int_{0 }^1 \,f_A(t,\bx,w) d w  - d_{23} \,f_R(t,\bx,\bv,u)  \int_{0 }^1  \,f_S (t,\bx,w) d w,
			\\[2mm]
			\mathcal N_C(\fvec) &\,=\,p_{42} \! \int_{\Gamma_{R}} \! \int_0^1 \! \int_0^1 \!\!
			\,f_A(t,\bx,u) \,f_R(t,\bx,\bv,w) du\,dw\,d\bv,
			\\[2mm]      
			\mathcal N_{E_1}(\fvec) &\,=\,- \widetilde{\mathcal N}_{E_1}(\fvec),
			\quad
			\mathcal N_{E_2}(\fvec) \,=\,\,  \widetilde{\mathcal N}_{E_1}(\fvec) - \bar{\mathcal N}_{E_2}(\fvec),\quad
			\mathcal N_{E}(\fvec) \,=\,\, \bar{\mathcal N}_{E_2}(\fvec),
		\end{aligned} 
		\ee
		with 
		\be
		\begin{aligned}\label{NEtil}
			\widetilde{\mathcal N}_{E_1}(\fvec) \,=\,&\,  b_{52}  \,f_{E_1}(t,\bx)   \int_{\Gamma_{R}}   \int_0^1 \,f_R(t,\bx,\bv,w)  dw\,d\bv,
			\\[2mm]
			\bar{\mathcal N}_{E_2} \,=\,&\,b_{62}  \,f_{E_2}(t,\bx)   \int_{\Gamma_{R}}   \int_0^1   \,f_R(t,\bx,\bv,w)  dw\,d\bv.
		\end{aligned}
		\ee
		
		\subsection{Operators corresponding to other processes}
		
		As in  \cite{della2022mathematical}, we include in our model the natural death of cell populations and decay of cytokines. We suppose that they occur at constant rate $d_I$, with $I\,=\,1,\,2,\,3,\,4$, corresponding to $A,\,R,\,S,\,C$, {respectively}. We also keep, as introduced in  \cite{della2022mathematical}, the constant input of self-antigen presenting cells. It represents, in fact, a constant increase of self-antigen presenting cells due to external environmental factors or life habits, represented by the parameter $\alpha$. Finally, we take into account the remyelination process performed by oligodendrocytes.
		Also, in this case, we look at the two phases of restoration. The first replacement of destroyed myelin with the partially sane
		\be \label{rem1} E+O\rightarrow E_2+O,\ee
		and the complete reformation of sane myelin
		\be \label{rem2} E_2+O\rightarrow E_1+O.\ee 
		The restoration rates for the two dynamics above are giv{en} by constants $r_I$, $I\,=\,5,6$ for \eqref{rem1} and \eqref{rem2}, respectively.
		{The} 
		operators tak{en} into account {in} these processes now read as 
		\be
		\begin{aligned}                
			\mathcal I_A(\fvec)&\,=\,\,\alpha - d_1  \,f_A(t,\bx,u),\qquad     
			\mathcal I_S(\fvec)\,=\,\,  - d_3  \,f_S(t,\bx,u) \\[2mm]
			\mathcal I_R(\fvec)&\,=\,\,- d_2  \,f_R(t,\bx,\bv,u),\qquad 
			\mathcal I_C(\fvec) \,=\,\, - d_{4} \,f_C(t,\bx,\bv)\\[2mm] 
			\mathcal I_{E_1}(\fvec) &\,=\,\, \widetilde{\mathcal I}_{E_2}(\fvec), \qquad 
			\mathcal I_{E_2}(\fvec) \,=\,\, -\widetilde{\mathcal I}_{E_2}(\fvec) + \bar{\mathcal I}_{E}(\fvec), \qquad {\mathcal I}_{E}(\fvec)\,=\,-\bar{\mathcal I}_{E}(\fvec),
		\end{aligned}  
		\ee
		with
		\begin{equation}
			\widetilde{\mathcal I}_{E_2}(\fvec) \,=\,\,  r_{5}  \,f_{E_2}(t,\bx) ,\qquad
			\label{IEtil} 
			\bar{\mathcal I}_{E} \,=\,\,r_{6}  \,f_{E}(t,\bx).   \end{equation}
		\subsection{Turning operators}\label{SubSecTurnOp}
		
		As pointed out previously, cytokines play a crucial role in recruiting self-reactive leukocytes in the central nervous system. Furthermore, they enhance the subsequent inflammation against cerebral tissues. 
		The chemical attraction operated by cytokines for immune cells, that detect the chemical concentration gradient and move in the direction
		{of} the source, has been widely shown by experiments\cite{taub1995chemotaxis}. 
		Several macroscopic models for chemotaxis, with similarities to the classical ones introduced in  \cite{keller1970initiation,patlak1953random}, can be found in literature, along with discussion of their analytical properties\cite{hillen2009user
		}. Some of them
		have been obtained from the kinetic {{theory modeling}}  by means of a turning operator, inspired by the ones proposed in  \cite{othmer2000diffusion,othmer2002diffusion}. Then, they have been extended in different ways, to reproduce the movement of cells through a velocity-jump process. Finally, these models have been applied to the description of a general autoimmune response in  \cite{immuneProc}. 
		Thus, we want to adapt the same procedure to the study of Multiple Sclerosis.
		In particular, inspired by Wang and Hillen in  \cite{wang2007classical}, we are going to adopt an analogous operator, with a more generalized form for the functions, including the volume-filling effect, expressed by a nonlinear squeezing probability (reproducing the elasticity of cells).To this aim, we start by setting a maximal value for the macroscopic density of self-reactive leukocytes, say $R_M>0$, and we {introduce} the functions $\varphi_0$ and $\varphi_1$ fulfilling the following properties.
		\begin{itemize}
			\item[(i)] Function $\varphi_1\in C^3[0,R_M]$ is such that $
			\varphi_1(0)\,=\,1,\,$ $
			\varphi_1(y) \in (0,1)\, \forall\, y\in (0,R_M)$ and
			$\varphi_1(y)\,=\,0\, \forall\, y \geq R_M$.
			\vspace{2mm}
			\item[(ii)] First and second derivatives of  the function $\varphi_1$ are such that $\varphi_1'(y)\leq 0$ with  $|\varphi_1'(y)|<+\infty$, $\varphi_1''(y)\leq 0$ for any $y\in [0,R_M]$.
			\vspace{2mm}
			\item[(iii)] Function $\varphi_0$ is related to $\varphi_1$ through the relation $\varphi_0(y)\,=\,\varphi_1(y)-\varphi_1'(y)\,y.$
		\end{itemize}

		{At this point, we write the turning operator for self-reactive leukocytes as
			\be 
			\mathcal {L}_R[\,f_C](f_R)(\bv) \,=\, \mathcal {L}_R^0(f_R)(\bv)+\mathcal {L}_R^1[\,f_C](f_R)(\bv),
			\label{eq:lr}
			\ee
			in order to differentiate the random  change in the velocity of the cells from the bias given by the chemotactic attraction of cytokines, as modeled in  \cite{othmer2002diffusion}, }
		More precisely, we take
		\be \nn
		\label{L0}
		\mathcal {L}_R^0(f_R)(\bv) \,=\,  \frac{\lambda}{\varphi_0(R)} \left( -\,f_R(\bv) + \int_{\Gamma_{R}}   T_R^0(\bv,\bvp)\,f_R(\bvp) d\bvp\right) .
		\ee
		The turning kernel $T_R^0$, as already introduced in  \cite{immuneProc}, gives the probability of a 
		{cell to change its velocity from $\bvp$ to $\bv$.}
		In this case, it reads as the uniform probability over the space of velocities (that here includes also the possibility of a change in both speed and direction, differently from origin{al} models as {in}\cite{othmer2000diffusion}).
		{T}hus, we have
		\be 
		\label{eq:T0}
		T_R^0\,=\,\frac{1}{\omega},
		\ee
		with $\omega\,=\,|\Gamma_R|$ {being} the measure of the velocity space. The parameter $\lambda>0$ represents the turning rate. 
		
		On the other hand, we define the operator $\mathcal {L}_R^1[\,f_C](f_R)$ as 
		\be \nn
		\label{L1} \mathcal {L}_R^1[\,f_C](f_R)(\bv) \,=\,\int_{\Gamma_{R}}\,\lambda\, T_R^1(\bv,\bvp,C)\,f_R(\bvp) d\bvp.
		\ee
		The kernel $T_R^1$ is {taken} in a more general form than the one used in  \cite{immuneProc}. Here, the subsequent velocity is a
		{byproduct} of the reorientation of the cell towards the gradient of the cytokines density, with a higher expectation of having a greater speed and with a dependence on the macroscopic density of self-reactive leukocytes. Then, we write
		\be
		\label{eq:T1}
		T_R^1(\bv,\bvp,C) \,=\, \gamma\,\varphi_1(R)\left(\hat{\bv}\cdot\hat{\bv}'\right)\left(\hat{\bv}'\cdot\nabla C\right)\,\frac{v}{\mathcal V}, \quad 
		\ee 
		having used the notation $\bv\,=\,v \hat{\bv}$, $|\hat{\bv}|\,=\,1$. The parameter
		$\gamma>0$ is the microscopic chemotaxis parameter.
		
		For the cytokines, instead, we only consider the random motion within the space of possible velocities, defining 
		\be 
		\label{LC}
		\mathcal {L}_C(f_C)(\bv) \,=\, \sigma \left( -\,f_R(\bv) + \int_{\Gamma_{C}}   T_C^0(\bv,\bvp)\,f_C(\bvp) d\bvp\right),
		\ee
		with $\sigma>0$. Again, we pick a uniform probability as kernel{, namely}
		\be \label{T0C}
		T_C^0(\bv,\bvp)\,=\,\frac{1}{\phi},\quad \phi\,=\,|\Gamma_C|.
		\ee
		We now provide the analytical properties of operators $\mathcal {L}_R^0(f_R)$, $\mathcal {L}_R^1[\,f_C](f_R)$  and $\mathcal {L}_C(f_C)$ that ensure the derivation of {a} reaction-diffusion macroscopic model.

		\begin{lemma}
			\label{lm:T} Let $\mathcal L_R[\,f_C]$ be the turning operator given {in} \eqref{eq:lr},  with  turning kernels $T_R^0$ and $T_R^1$ defined in \eqref{eq:T0} and \eqref{eq:T1}, respectively, and let $\mathcal L_C$ be the turning operator defined in \eqref{LC} having as turning kernel $T_C^0$ given in \eqref{T0C}.
			\begin{itemize}
				\item[(i)] Being the turning kernels $T_R^0$ and  $T_C^0$  constant and positive, and holding (for I\,=\,R,\,C)
				\be \nn
				\int_{\Gamma_I} \!\! T_I^0(\bv,\!\bvp) d\bv
				\,=\,\dis\int_{\Gamma_I} \!\! T_I^0(\bv,\!\bvp) d\bvp
				\,=\,\dis\int_{\Gamma_I} \int_{\Gamma_I} \!\! (T_I^0)^2(\bv,\!\bvp) d\bv d\bvp \,=\, 1 
				\label{eq:PT0}
				\ee
				{then} the equation $\mathcal L_I^0(f) \!\,=\,\! g$ {has} a unique solution $f\in L^2\left(\Gamma_I\right)$,
				owning the solvability condition 
				\be\label{SolvCond1}
				\int_{\Gamma_I} f \,d\bv \,=\, 0 \Leftrightarrow \int_{\Gamma_I} g \,d \bv \,=\, 0 .
				\ee
				Moreover, the inverse operator of  $\mathcal L_I^0$  {corresponds to} the multiplication  by $-{\varphi_0(R)}/{\lambda}$ for $\mathcal L_R^0$ and by $-{1}/{\phi}$ for $\mathcal L_C^0$.
				\vspace{2mm}
				\item[(ii)] The kernel  $T_1(\cdot,\cdot,C)$  is a $L^2${-fu}nction of the  space $\Gamma_R\times \Gamma_R $ for any density $C$ of cytokines and we have
				\be \nn
				\int_{\Gamma_R} T_R^1(\bv,\bvp,C) d\bv\,=\,0 .
				\label{eq:PT1}
				\ee
			\end{itemize} 
		\end{lemma}
		{The spectral properties listed above constitute a particular case of more general results provided in  \cite{othmer2000diffusion,othmer2002diffusion}, to which we address the reader for the proofs.}

		\section{Diffusive limit}
		\label{Sec4}
		
		{In this section, we aim to apply common asymptotic methods commonly implemented in the kinetic {{theory}} framework to obtain a diffusive limit {of the system \eqref{eq:kin1}-\eqref{eq:kin2}}.}

		\subsection{Time scaling}
		
		The key assumption is to consider the different time scales at which the various processes occur. For this reason, we take a small characteristic parameter $\varepsilon$ and we claim the following
		{assertions.}
		\begin{itemize}
			\item[(i)] Velocity-jump processes are the dominant process in the dynamics, namely of order $\varepsilon^{-1}$.
			\item[(ii)] The velocity-jump process induced by cytokines on self-reactive leukocytes occurs at a slower time scale (of order $\varepsilon^1$) with respect to the random changes in velocity (order $\varepsilon^0$).
			\item[(iii)]  Conservative and non-conservative interactions, constant input of self-antigen presenting cells, and natural death of cells are of order $\varepsilon$.
			\item[(iv)]  Processes involving initial consumption of myelin \eqref{dem1} and total restoration \eqref{rem2} are slower (of order $\varepsilon^2$), than the complete consuming of myelin \eqref{dem2} and the initial restoration \eqref{rem1}.
		\end{itemize}
		
		Assumptions (i) and (ii) are inspired by the model proposed in  \cite{othmer2002diffusion}, while (iii) is an application of the general setting implemented by authors in  \cite{bellomo2016multiscale}. Finally, {the} time scaling proposed in (iv) for the myelin consumption comes from heuristic considerations (based on observations reported in  \cite{zhao1996mri}) and from the choice of focusing on the active stage of lesions.  For this last reason, we take the time scale $o(\varepsilon^1)$ and proceed analogously to derivations previously conducted in the framework of kinetic theory for gas mixtures\cite{bisi2006reactive}. Specifically, given that we are focusing on the time scale related to slow processes listed in assumptions (iii), we apply the same scaling in front of the temporal derivative.
		Therefore, from \eqref{eq:kin1}-\eqref{eq:kin2}, we obtain the following rescaled system
		
		\begin{equation}
			\begin{cases}
				\dis \varepsilon\,  \,\frac{\pa \,f_A}{\pa t} 
				\,=\, \varepsilon\,\,\mathcal  N_A(\fvec)+ \varepsilon\,\,\mathcal I_A(\fvec) ,
				\label{BolA}
				\\[3mm]
				\dis \varepsilon\,    \,\frac{\pa \,f_S}{\pa t}
				\,=\, \varepsilon\,\,\mathcal  N_S(\fvec) + \varepsilon\,\,\mathcal I_S(\fvec) ,
				\\[3mm]
				\dis \varepsilon\,  \,\frac{\pa \,f_R}{\pa t} + \bv\cdot\nabla_{\bx}\, \,f_R
				\,=\, \frac1\varepsilon\, \mathcal {L}^{\varepsilon\,}_R[\,f_C](f_R) + \varepsilon\, \mathcal N_R(\fvec) +\varepsilon\, \mathcal I_R(\fvec) ,
				\\[3mm]
				\dis \varepsilon\,  \, \frac{\pa \,f_C}{\pa t}+\bv\cdot\nabla_{\bx}\, \,f_C
				\,=\, \frac1\varepsilon\, \mathcal {L}_C(f_C) + \varepsilon\, \, \mathcal { N_C}(\fvec) + \varepsilon\,\,\mathcal I_C(\fvec),
				\\[3mm]
				\dis \varepsilon\,  \, \frac{\pa \,f_{E_1}}{\pa t}
				\,=\, -\varepsilon^2\, \, \widetilde{\mathcal N}_{E_1}(\fvec) + \varepsilon^2\,\,\widetilde{\mathcal I}_{E_2}(\fvec),
				\\[3mm]
				\dis \varepsilon\,  \, \frac{\pa \,f_{E_2}}{\pa t}
				\,=\, \varepsilon^2\, \, \widetilde{\mathcal N}_{E_1}(\fvec) - \varepsilon^2\,\,\widetilde{\mathcal I}_{E_2}(\fvec)-\varepsilon\, \, \bar{\mathcal N}_{E_2}(\fvec) + \varepsilon\,\,\bar{\mathcal I}_E(\fvec),
			\end{cases}
		\end{equation}
		with
		\be \nn
		\mathcal {L}^{\varepsilon}_R[\,f_C](f_R)(\bv) \,=\, \mathcal {L}_R^0(f_R)(\bv)+\varepsilon\,\mathcal {L}_R^1[\,f_C](f_R)(\bv) .
		\label{eq:lrSca}
		\ee
		
		{
			We have omitted the equation for the destroyed myelin $E$ since it is easy to check that the total quantity of myelin (sane, partially sane, and destroyed) is constant in time. {We consider also the initial total density constant in space, consequently $E_1(t,\bx)+E_2(t,\bx)+E(t,\bx)\,=\,\bar E$.}}
		We also emphasize that the different powers of $\epsilon$ preceding each term in the above system indicate their order of dominance concerning the considered time scale.
		
		At this point, we aim to obtain a macroscopic description of the model, with quantities depending only on variables $t$ and $\bx$. 
		We start by performing a Hilbert expansion in $\varepsilon$ of the distribution functions,
		obtaining
		\be
		\begin{aligned} 
			\,f_A(t,\bx,u) \,=\,& { \,f_A^0(t,\bx,u)} + \varepsilon \,f_A^1(t,\bx,u) + \varepsilon^2\,f_A^2(t,\bx,u) + O(\varepsilon^3),
			\label{expfA} 
			\\[2mm]
			\,f_S(t,\bx,u) \,=\,& { \,f_S^0(t,\bx,u)} + \varepsilon \,f_S^1(t,\bx,u) + \varepsilon^2\,f_S^2(t,\bx,u) + O(\varepsilon^3),
			\\[2mm]
			\,f_R(t,\bx,\bv,u) \,=\,& { \,f_R^0(t,\bx,\bv,u)} + \varepsilon \,f_R^1(t,\bx,\bv,u) + \varepsilon^2\,f_R^2(t,\bx,\bv,u) + O(\varepsilon^3),
			\\[2mm]
			\,f_C(t,\bx,\bv) \,=\,& { \,f_C^0(t,\bx,\bv)} + \varepsilon \,f_C^1(t,\bx,\bv) + \varepsilon^2\,f_C^2(t,\bx,\bv) + O(\varepsilon^3),
			\\[2mm]
			\,f_{E_1}(t,\bx) \,=\,& { \,f_{E_1}^0(t,\bx)} + \varepsilon \,f_{E_1}^1(t,\bx) + \varepsilon^2\,f_{E_1}^2(t,\bx) + O(\varepsilon^3),
			\\[2mm]
			\,f_{E_2}(t,\bx) \,=\,& { \,f_{E_2}^0(t,\bx)} + \varepsilon \,f_{E_2}^1(t,\bx) + \varepsilon^2\,f_{E_2}^2(t,\bx) + O(\varepsilon^3).
		\end{aligned}
		\ee
		Observing system \eqref{BolA}, it is noticeable that an expansion up to the second order is redundant for certain distribution functions. However, we retain the notation above for the sake of much generality.
		{Similarly to what has been assumed in previous works\cite{othmer2000diffusion}, we consider} 
		that the {whole} mass is concentrated in the $0$-th order terms, {i.e.}, for $k\geq1$ 
		\be\nn
		\int_0^1 \,f_I^k(t, \bx, u)\,du\,=\,0, \ \  I\,=\,A,\,R,\,S,\ee 
		\be\nn
		\int_0^1 \,f_R^k(t, \bx,\bv,u)\,du\,=\,\int_0^1 \,f_R^k(t, \bx,\bv,u)\,d\bv\,=\,0,
		\ee
		\be\nn
		\int_0^1 \,f_C^k(t, \bx,\bv)d\bv\,=\,0, \quad \,f_{E_1}^k(t, \bx)\,=\,0, \quad \,f_{E_2}^k(t, \bx),=\,0.
		\ee

		\subsection{Diffusion equation for self-reactive leukocytes}
		\label{ssec:r}
		
		For convenience, we {start by} considering the third equation in {\eqref{BolA}} for self-reactive leukocytes, and insert  the corresponding expansions given in
		\eqref{expfA}. 
		{Then, w}e equal {the} terms of order $\varepsilon^0$, $\varepsilon^1$ and $\varepsilon^2$, 
		obtaining the {following} equations
		\begin{align}
			\mathcal {L}_R^0( \,f_R^{0})& \,=\, 0,
			\label{eq:e0}
			\\[2mm]
			\bv\cdot\nabla_{\bx}\, \,f_R^0 &\,=\, \mathcal {L}_R^0( \,f_R^{1})+\mathcal {L}_R^1[\,f_C]( \,f_R^{0}),
			\label{eq:e1}
			\\[2mm]
			\frac{\pa \,f_R^0}{\pa t} \!+\! \bv\!\cdot\!\nabla_{\bx} \,f_R^1 	&\,=\,  \mathcal {L}_R^0( \,f_R^{2})
			+ \mathcal {L}_R^1[\,f_C]( \,f_R^{1})
			+ \mathcal N_R(f_A^0,\,f_R^0,\,f_S^0) + \mathcal I_R(f_R^0).
			\label{eq:e2}
		\end{align}
		Thanks to the spectral properties of operator $\mathcal {L}_R^0(f_R)(\bv)$ stated in Lemma \ref{lm:T}, 
		we can {determine the} explicit {expressions} of {the} order  $\varepsilon^0$ and $\varepsilon^1$ terms in the expansion of $R$ in \eqref{expfA}, 
		that are
		\be \nn
		\,f_R^0{ (t,\bx,\bv,u)} \,=\, \rho_R(t,\bx,u) ,
		\ee
		and 
		\be \nn
		\,f_R^{1} { (t,\bx,\bv,u)} 
		\,=\, - \frac{\varphi_0(R)}{\lambda} \bv \cdot \nabla_{\bx} \rho_R
		+  \rho_R\int_{\Gamma_{R}} T_R^1(\bv,\bvp,C) d\bvp.
		\ee
		{Then, s}ubstituting the terms {$\,f_R^0$ and $\,f_R^1$} in equation \eqref{eq:e2}, 
		in order to recover the term $ \,f_R^2$, we have  {first} to apply the solvability condition \eqref{SolvCond1} that leads to
		\begin{align}\nn 
			&\omega \frac{\pa \rho_R}{\pa t}-
			\nabla_{\bx}\cdot\left(\int_{\Gamma_{R}}\frac{\varphi_0(R)}{\lambda} \bv \otimes \bv\right)\cdot \nabla_{\bx}\, \rho_R  +\nabla_{\bx}\cdot\left[\rho_R\int_{\Gamma_{R}}\int_{\Gamma_{R}} \bv \, T_R^1(\bv,\bvp,C) d\bvp\,d\bv\right]\\[2mm]
			&\qquad\,=\,
			\left[\mathcal N_R(f_A^0,\rho_R,\,f_S^0) \!+\! \mathcal I_R(\rho_R)\right]\omega.
			\label{eq:ror}
		\end{align}
		{A}fter comput{ing the} integrals, 
		{equation \eqref{eq:ror}}
		becomes the reaction-diffusion equ{a}tion for $\rho_R$,
		\begin{align}
			&\nn \frac{\pa \rho_R}{\pa t} - \nabla_{\bx}\,\cdot  \left[D_R\,\varphi_0(R)\,\nabla_{\bx}\,\,\rho_R -\chi\,\varphi_1(R)\,\rho_R \,\nabla_{\bx}\, \,C\right]
			\\[2mm]
			& \qquad \,=\, \mathcal N_R(f_A^0,\rho_R,\,f_S^0) \!+\! \mathcal I_R(\rho_R),
			\label{eq:rorr}
		\end{align}
		with the diffusion coefficient $D_R$ and the chemotactic parameter $\chi$ {given by}
		\be \label{DRChi}
		D_R \,=\, \frac{\mathcal V^2}{(n+2)\lambda} \quad \mbox{ and } \quad 
		\chi \,=\, \gamma\,\omega\,\frac{ \mathcal V}{(n+1)^2},
		\ee
		respectively. 
		
		\subsection{Evolution equations for self-antigen presenting cells and immunosuppressive cells}

		{L}et us consider equations for the populations of self-antigen presenting cells and immunosuppressive cells in \eqref{BolA},
		and insert the corresponding expansions given in \eqref{expfA}.
		Equaling the same order terms {in $\varepsilon$}, we get only the equation{s} involving the $\varepsilon^0$ terms, 
		that are
		\begin{align}
			\dis  \,\frac{\pa \,f_A^0}{\pa t} &
			\,=\, \mathcal  N_A[\,f_A^0,\rho_R,\,f_S^0]+ \mathcal I_A(f_A^0) ,
			\label{BolAze}
			\\
			\dis  \,\frac{\pa \,f_S^0}{\pa t} &
			\,=\, \mathcal  N_S[\,f_A^0,\rho_R,\,f_S^0]\,+\, \mathcal I_S(f_S^0).
			\label{BolSze}
		\end{align}
		{ Such equations describe the behavior in time of the distribution functions of self-antigen presenting cells and immunosuppressive cells due to proliferative and destructive processes.
		}	
		
		\subsection{Diffusion equation for the cytokines}
		
		The procedure applied {in Subsection \ref{ssec:r}} to the equation for self-reactive leukocytes can be analogously followed to recover the macroscopic equation for cytokines density. {Then, w}e insert {expansions} for $\,f_R$ and for $\,f_C$ in  \eqref{expfA}
		{in the  equation} for cytokines in \eqref{BolA}, obtaining also, in this case, an equality for each order of $\varepsilon$, 
		{i.e.}, 
		\begin{align}
			\mathcal {L}_C^0( \,f_C^{0}) &\,=\, 0,
			\label{eq:e0C}
			\\[2mm]
			\bv\cdot\nabla_{\bx}\, \,f_C^0 &\,=\, \mathcal {L}_C^0( \,f_C^{1}),
			\label{eq:e1C} 
			\\[2mm]
			{\frac{\pa \,f_C^0}{\pa t} \!+\! \bv\!\cdot\!\nabla_{\bx} \,f_C^1 }&
			\,=\,  \mathcal {L}_C^0( \,f_C^{2}) \!+\! \mathcal N_C(f_A^0,\rho_R,\,f_C^0) \!+\! \mathcal I_C(f_C^0).
			\label{eq:e2C}     
		\end{align} 
		Holding again Lemma \ref{lm:T}, we find the first two terms of the expansion, 
		that turn to be
		\be\nn
		\,f_C^0(t,\bx,\bv)\,=\,C(t,\bx),\quad \,f_C^1(t,\bx,\bv)\,=\,-\frac{1}{\sigma}\bv\cdot\nabla_{\bx}C(t,\bx),
		\ee
		so that the solvability condition \eqref{SolvCond1}, applied to {equation} \eqref{eq:e2C},  gives
		\begin{align}
			\frac{\pa C}{\pa t}-&
			D_C\,\Delta_{\bx}\,\,C  \,=\,
			\mathcal N_C[\,f_A^0,\rho_R,\,C] \,+\, \mathcal I_C(C),
			\label{eq:C}
		\end{align}
		with {the diffusion} coefficient for cytokines {given by}
		$ \displaystyle
		D_C \,=\, \frac{\mathcal W^2}{(n+2)\sigma}$.

		\subsection{Evolution equations for sane and partially sane myelin}
		
		Finally, we consider equations for the sane and partially sane myelin. Inserting expansions for $f_R$, $f_{E_1}$ and equation in the penultimate equation of \eqref{BolA} and equaling same order terms, we find 
		\be
		\dis  \, \frac{\pa \,f_{E_1}}{\pa t}
		\,=\, 0,\quad
		- \, \widetilde{\mathcal N}_{E_1}[\rho_R,\,f_{E_1}] + \,\widetilde{\mathcal I}_{E_2}[\,f_{E_2}]\,=\,0,\label{E1cost}
		\ee
		and, recalling expressions for $\widetilde{\mathcal N}_{E_1}$ and $\widetilde{\mathcal I}_{E_2}(\fvec)$ given in \eqref{NEtil} and \eqref{IEtil}, respectively, we have the relation 
		\be r_{5} \,f_{E_2}(t,\bx) - b_{52}\,f_{E_1}(t,\bx) R(t,\bx)\,=\,0.\label{relE1E2}\ee
		On the other hand, as we insert expansions in the very last equation of  \eqref{BolA} we get at leading order
		\be
		\dis  \, \frac{\pa \,f_{E_2}}{\pa t}\,=\,
		\, -\bar{\mathcal N}_{E_2}[\rho_R,\,f_{E_2}] + \,\bar{\mathcal I}_E[\,f_{E_2}].\label{eq:E2}
		\ee
		
		Macroscopic equations for the constituents of the model are recovered integrating equations \eqref{eq:rorr}, \eqref{BolAze}, \eqref{BolSze} with respect to the activity variable $u$ and writing all operators explicitly in \eqref{eq:rorr}, \eqref{BolAze}, \eqref{BolSze}, \eqref{eq:C} and \eqref{eq:E2} (using the notation $E_I\,=\,\,f_{E_I}$). 
		Moreover, since we are interested in the portion of destroyed myelin, we exploit the relation $E(t,\bx)\,=\,\bar E-E_1(t,\bx)-E_2(t,\bx)$ from which, using \eqref{relE1E2}, we can write
		\be\nn
		E_2(t,\bx)\,=\,(\bar E - E(t,\bx))\frac{b_{52}\,R(t,\bx)}{r_5+b_{52}\,R(t,\bx)}.
		\ee
		In addition, \eqref{E1cost} provides ${\pa_t E}\,=\,-{\pa_t E_2}$.

		\subsection{Limiting equations}
		
		{As a result of the previous subsections, in the formal limit $\varepsilon \rightarrow 0$, we obtain the following system of equations}
		\be
		\begin{cases}
			\label{eq:macA}
			\dfrac{\pa \,A(t,\bx)}{\pa t} \,=\, \alpha + p_{12} \,A(t,\bx) \,R(t,\bx) - d_{13} \,A(t,\bx) \,S(t,\bx) - d_1 \,A(t,\bx), 
			\\[3mm]
			\dfrac{\pa \,S(t,\bx)}{\pa t} \,=\, p_{31} \,S(t,\bx) \,A(t,\bx) - d_3 \,S(t,\bx),
			\\[3mm]
			\dfrac{\pa \,R(t,\bx)}{\pa t} \,=\,
			\nabla_{\bx}\,\cdot\left[D_R\,\varphi_0(\,R(t,\bx))\nabla_{\bx}\,\,\,R(t,\bx) - \chi \,\varphi_1(\,R(t,\bx))\,\,R(t,\bx)  \,\nabla_{\bx}\,C\right]\\ 
			\qquad\qquad\qquad+ p_{21} \,R(t,\bx) \,A(t,\bx) - d_{23} \,R(t,\bx) \,S(t,\bx) - d_2 \,R(t,\bx),
			\\[3mm]
			\dfrac{\pa C(t,\bx)}{\pa t}\,=\,D_C\, \Delta_{\bx}\,C(t,\bx) +p_{C2}\,A(\bx,t)\,R(\bx,t)- d_C C(\bx,t),
			\\[3mm]
			\dfrac{\pa E(t,\bx)}{\pa t}\,=\,(\bar E - E(t,\bx))\dfrac{b_{52}\,b_{62}\,R(t,\bx)}{r_5+b_{52}\,R(t,\bx)}\,R(t,\bx)-r_{6}E(t,\bx).
		\end{cases}
		\ee
		The system obtained {above} describes the global behavior of the populations considered in our model. 
		{Its linear stability will be studied}
		in the following sections, allowing us to understand the interplay of different agents. In particular, the motion of self-reactive leukocytes, driven by chemotactic attraction, may lead to a patterned-in-space consumption of myelin, as described by the last equation. To this aim, we shall consider the problem with  non-negative initial data 
		\begin{equation}\label{InCond}\nn
			\mathbf{U}(0,\bx)\,=\,\mathbf{U}_0(\bx)\geq 0, \mbox{ with } \mathbf{U}(t, \bx) \,=\, (A, S, R, C, E), \mbox{ and } R_0(\bx)\leq R_M,
		\end{equation}
		and no-flux at the boundary, 
		\be\label{NoFlux}\nn
		\Big(\varphi_0(R)\,\nabla_{\bx}\, R-{\chi} \,\varphi_1(R) R\,\nabla_{\bx}\,C \Big)\cdot {\bf \widehat n}\,=\,0,\quad \nabla_{\bx} C\cdot {\bf \widehat n}\,=\,0,
		\ee
		with  ${\bf \widehat n}$ {being} the external unit normal to {the boundary} $\pa\Gamma_{\bx}$.
		
		Before proceeding, we aim to provide insights into the characteristics of the solution, focusing particularly on its boundedness and positivity.

		\begin{remark}
			{ We consider, at first, the system \eqref{eq:macA} without the diffusive and chemotactic terms. In this case, it is possible to extend results obtained in  \cite{della2022mathematical} and ensure boundedness and positivity of the solution in the space-dependent case, through a time-step discretization. We perform this calculation in \ref{App} of the present paper. }Considering the full model case, we can not apply the same reasoning due to the particle flux created by chemotaxis and diffusion. Anyway, since chemotaxis and diffusion do not occur for $A$, $S$, and $E$, we can show the boundedness and positivity of the solution for a reduced version of the system. Such a reduction is analogous to the one proposed in  \cite{smith2018model}. It can be obtained by writing the non-diffusing species ($A$, $S$, $E$) in terms of the diffusing ones ($R$, $C$), by imposing 
			$\partial_t A\,=\,\partial_t S\,=\,\partial_t E\,=\,0$ and obtaining
			\begin{align}\nn
				A&\,=\,\frac{d_3}{p_{31}} ,\\\nn
				S&\,=\, \frac{\alpha p_{31}-d_1\, d_3}{d_3\, d_{13}} + \frac{p_{12}}{d_{13}}R ,\\\nn
				E&\,=\,\frac{R^2\,\bar E\,b_{52}\,b_{62}}{r_6\,(r_5+R\,b_{52})+R^2\,b_{52}\,b_{62}}
				.
			\end{align}
			Replacing these expressions in equations for $R$ and $C$ in \eqref{eq:macA}, we get 
			\be
			\begin{cases}\nn
				\dfrac{\pa R(t,\bx)}{\pa t} &\,=\,
				\nabla_{\bx}\,\cdot\left[D_R\,\varphi_0(R(t,\bx))\nabla_{\bx}\,\,R(t,\bx) - \chi \,\varphi_1(R(t,\bx))\,R(t,\bx)  \,\nabla_{\bx}\,C\right]+ \,f(R),
				\\[3mm]\nn
				\dfrac{\pa C(t,\bx)}{\pa t}&\,=\,D_C\, \Delta_{\bx}\,C(t,\bx) +g(R,\,C) ,
			\end{cases}
			\ee
			where the functions $\,f(R)$ and $g(R,\,C)$ are defined by
			\begin{align}\nn
				\,f(R) \,=\, aR\left(1-\frac{b}{a}R\right) , \quad
				g(R,\,C) \,=\, \frac{p_{C2}\,d_3}{p_{31}}R(\bx,t)- d_C C(\bx,t) ,
			\end{align}
			with
			{\begin{align}\nn
					a\,=\, \frac{p_{21}\,d_2}{p_{31}} - d_{23}\left(\frac{\alpha\,p_{31} - d_1\,d_3}{d_2\,d_{13}}\right) - d_2,\quad
					b\,=\, \frac{p_{12}\,d_2}{d_{13}} \ .
			\end{align}}
			In this regime, the system $(R,\, C)$ is analogous to the one presented in  \cite{wang2007classical} that, as shown in that paper, respects the conditions for positivity, boundedness, and uniqueness as long as $a>0$. This regime includes solutions close to equilibria of the non-diffusive and non-chemotactic model, that we show to be positive and bounded in \ref{App}, and are relevant to the next section. { For a more complete overview of qualitative analysis for global classical solvability, boundedness, and large time behavior in diverse types of chemotaxis syste{ms,} we address the reader to  \cite{ke2022analysis}; an enriched review of models including interaction with a time-evolving external system is given in  \cite{bellomo2022chemotaxis}. }
		\end{remark}

		\section{Turing instability analysis}
		\label{Sec5}
		
		In this section, we are going to show how the macroscopic system \eqref{eq:macA}, under a particular choice of parameters, exhibits spatial patterns representing the {appearance} (and possible reconstruction) of myelin plaques. 
		
		\subsection{Adimensionalization}
		We  first adimensionalyze the system, adopting the change of variables
		\be\nn \widetilde t=d_3 t,\qquad \widetilde{\bx}=\sqrt{\frac{d_3}{Dr}}\bx.\ee 
		Then, we introduce the adimensionalized quantities as follows
		\be \nn
		\widetilde A=\frac{A d_3}{\alpha},\quad
		\widetilde S=\frac{S d_{13}}{d_3},\quad
		\widetilde R=\frac{R}{R_M},\quad
		\widetilde C=\frac{C d_3^2}{{p_{C2}}\alpha\,R_M},\quad
		\widetilde E=\frac{E}{\bar E}.
		\label{eq:ndimdens}
		\ee
		Defining the new coefficients of the model {as}
		\be 
		\nn
		\beta=\frac{R_M\,p_{12}}{d_3},\quad
		\zeta=\frac{d_1}{d_3},\quad
		\mu =\frac{p_{31}\alpha}{d_3^2},\quad
		\delta=\frac{D_C}{D_R} ,\quad
		\tau=\frac{d_C}{d_3},\quad
		\xi=\chi\,\frac{ p_{C2}\alpha\, R_M}{D_R d_3^2\,},
		\ee
		\be
		\label{CoefAdim}
		\eta=\frac{p_{21}\alpha}{d_3^2},\quad
		\phi=\frac{d_{23}}{d_{13}},\quad 
		\theta=\frac{d_2}{d_3},\quad
		\Theta=\frac{b_{62}}{d_3},\,\quad
		\Omega=\frac{r_5}{R_M\,b_{52}},\quad
		\Xi=\frac{r_6}{d_3},\quad
		\ee
		and functions 
		\be\nn
		\Phi_k(x)\,=\,\varphi_k(x\,R_M),\quad k\,=\,0,1,
		\ee
		we get the {following dimensionless equations from} system
		{ \eqref{eq:macA}},
		\be
		\begin{cases}
			\label{eq:rdsA} \dfrac{\pa A}{\pa t}\,=\,1+\beta\,A\,R - A\,S - \zeta A , \\[3mm]
			\dfrac{\pa S}{\pa t}\,=\,\mu\,A\,S - \,S , \\[3mm]
			\dfrac{\pa R}{\pa t}\,=\,\nabla_{\bx}\cdot\left(\Phi_0(R)\,\nabla_{\bx}\, R-{\xi} \,\Phi_1(R) R\,\nabla_{\bx}\,C\right)+\eta\,A\,R - \phi\,R\,S - \theta\,R , \\[3mm]
			\dfrac{\pa C}{\pa t}\,=\,\delta \Delta_{\bx}\,C+A\, R-\tau\,C,\\[3mm]
			\dfrac{\pa E}{\pa t}\,=\, \dfrac{\Theta\,R}{\Omega+R}\,R\,(1-E)-\Xi\,E ,
		\end{cases}
		\ee
		{where we have renamed} the tilde-labeled densities {by} removing the tilde.

		We study the system { \eqref{eq:rdsA} by} {considering} initial data 
		\begin{equation}\label{InCond2}\nn
			\mathbf{U}(0,\bx)\,=\,\mathbf{U}_0(\bx)\geq 0, \mbox{ with } R_0(\bx)\leq 1,
		\end{equation}
		and imposing zer{o}-flux at the boundary,
		{i.e.}  
		\be\label{ZeroFlux}
		\Big(\Phi_0(R)\,\nabla_{\bx}\, R-{\xi} \,\Phi_1(R) R\,\nabla_{\bx}\,C \Big)\cdot {\bf \widehat n}\,=\,0,\quad \nabla_{\bx} C\cdot {\bf \widehat n}\,=\,0.
		\ee
		The formation of patterns due to Turing instability\cite{turing52} occurs when a spatially homogeneous steady state becomes unstable 
		{}due to the addition of  
		the diffusive terms. The stability properties of {the sub}system composed by the first four equations of \eqref{eq:rdsA}, without the diffusion, have already been studied for a general case of autoimmunity in  \cite{immuneProc}. Thus, we are going to show that they also hold in the case in which the model is applied to the particular case of Multiple Sclerosis. 
		
		First of all, we {can specify} a biologically relevant equilibrium for the system \eqref{eq:rdsA}	in spatially homogeneous conditions. Then we find
		$U_1\,=\, (A_ 1, \,S_1,\,R_1,\,C_1,\,E_1)$ that writes
		\be\label{U}
		U_1 \,=\, \,\left(\frac{1}{\,\mu}, 
		\,\frac{ \eta -\theta \,\mu}{\mu \,\phi},\,
		\frac{-\theta \,\mu\,  + \eta\,+ \mu\,\phi\,( \zeta\,-\mu  ) }{\beta\, \mu\, \, \phi},\,          
		\frac{1}{\mu\tau}\,R_1,\,\frac{R_1^2 \Theta}{R_1^2\, \Theta + 
			R_1\, \Xi + \Xi\, \Omega} \right),
		\ee
		belonging to the set
		\be \nn
		\mathcal E\,=\,\Big\{ A(t,\bx)>0,\, S(t,\bx)> 0,\, 0< R(t,\bx)\leq 1,\, C(t,\bx)> 0,\, E(t,\bx)> 0\Big\}.
		\ee
		
		We point out that equilibrium $U_1$ defined in \eqref{U} actually belongs to $\mathcal E$ provided {that} the following conditions on the parameter $\theta$ {hold,}
		\be\label{EqPos}
		\theta < \frac{\eta}{\mu}  
		\quad\mbox{and}\quad
		\bar\theta > \theta > \bar\theta-\beta\phi,
		\quad\mbox{with}\quad
		\bar\theta:\,=\,\frac{\eta}{\mu}+\phi(\zeta-\mu).
		\ee
		For convenience, we will set {the} parameters in such a way {that} $\zeta<\mu$, so that the second condition in \eqref{EqPos} implies the first one.
		
		The following step is {performed in order} to linearize the system \eqref{eq:rdsA} around equilibrium $U_1$, getting 
		\be \nn
		\frac{\pa \bf W}{\pa t} \,=\, {\mathbb A} {\bf W},
		\quad \mbox{for} \quad
		{\bf W}\,=\,\left(\begin{array}{c}A- A_1,\,S-S_1,\,R- R_1,\,C-C_1,\,E-E_1\end{array}\right)^\intercal ,
		\label{eq:lins}
		\ee 
		{where t}he Jacobian matrix $\mathbb A$ {is given by}
		\be \nn
		{\mathbb A}\,=\,\left(
		\begin{array}{ccccc}
			&&& 0 & 0\\
			&{\mathbb J}& & 0 &0 \\
			&&&0&0\\
			0&0&1&-\tau &0\\
			0&0&
			\dfrac{R_1\,\Theta\, \Xi\, (R_1 +2 \,\Omega)}{(R_1 +
				\Omega) (R_1\,(R_1\,\Theta + \Xi)+\Xi\,\Omega)}& 
			0 & -\Xi - \dfrac{R_1^2\,\Theta}{R_1+\Omega}
		\end{array}
		\right),
		\ee
		\mbox{with}
		\be \nn
		{\mathbb J}\,=\,
		\left(
		\begin{array}{ccc}
			-{ \mu} & -\dfrac{1 }{\mu}& \dfrac{\beta }{\mu} \\ &&\\
			\dfrac{ \eta -\theta \,\mu }{\phi} & 0 & 0 \\&&\\
			R_1\,\eta   & -R_1\,\phi & 0
		\end{array}
		\right).
		\ee
		{The e}igenvalues of the matrix $\mathbb A$ are $-\tau$ and $-\Xi - { \Theta}/({\beta\,\phi +\beta^2\,\phi^2\, \Omega})$, 
		along with the eigenvalues of {matrix} $\mathbb J$, having as characteristic polynomial $$
		P_{{\mathbb J}}(\lambda) \,=\,  -\left(\lambda^3 +A_1 \lambda^2+
		A_2\,\lambda\,+A_3\right),$$ with coefficients
		\be\nn
		A_1\,=\,  {\mu},\quad A_2\,=\,\left(-\frac{R_1 \,\beta \,\eta}{\mu} - \frac{\theta}{\phi} + \frac{\eta}{\mu\, \phi}\right),\quad A_3\,=\,
		R_1\,\beta\,\frac{\eta-\theta\,\mu}{\mu}.\ee
		Accordingly to the Routh-Hourwitz criterion\cite{gantmacher1964theory}, {the} eigenvalues {of ${\mathbb J}$} have all negative real part if and only if $A_1>0$ (automatically satisfied), $A_3>0$ (holds true thanks to \eqref{EqPos}) and $A_1\,A_2>A_3$. 
		This last requirement reads as
		\be\nn
		\frac{\mu(\eta-\theta\,\mu)-R_1 \,\beta (\theta\,\,\mu \,\nu + \eta (\mu + \nu)) \,\phi}{\mu\,\nu\,\phi}>0,
		\ee
		that, by \eqref{U}, becomes
		\be \label{ThetaPM}
		-\theta^2+\theta\,\left(-\eta - 2\frac{ \eta}{\mu} + \mu +\phi\,(\zeta - \mu)\right)+\eta -\frac{\eta^2}{\mu^2}+\frac{ \eta\,\phi (1 + \mu) (\mu-\zeta) -\eta^2}{\mu}>0.
		\ee
		It means that, under suitable choice of param{et}ers, it is possible to individuate an interval $(\theta_-,\theta_+)$, 
		with 
		\begin{align}\nn
			\theta_{\pm}\,=\,&\frac{\eta}{\mu} + \frac12\,\left(-\mu\,(1 + \phi) + (\eta + \zeta\,\phi)\right)\\[2mm]\nn
			&\quad \pm\frac12\sqrt{\eta^2 + 2 \eta\,\mu (\phi-1) - 
				2 \zeta\,\eta\,\phi + (\mu - \zeta\,\phi + \mu\,\phi)^2},
		\end{align}
		such that {the} equilibrium $U_1$ is unstable in spatially homogeneous conditions for any choice of 
		$\theta$ satisfying
		\be\nn
		\min(\bar\theta-\beta\,\phi, \theta_-)<\theta<\max(\bar\theta, \theta_+).
		\ee
		Moreover, we observe that for {$\theta\,=\,\theta_+,\, \theta_-$} we have $A_1>0$, $A_2>0$ and $A_1\,A_2\,=\,A_3$, which implies that a Hopf bifurcation\cite{marsden1976hopf} occurs and the solutions of the spatially homogeneous system {becomes} periodic in time with period
		\be\nn
		T\,=\,\frac{2\pi}{\sqrt{A_2}}.
		\ee
		
		{In view of investigating Turing instability,} we consider the complete system with the diffusion {effects, namely}
		\eqref{eq:rdsA}, {and linearize the system} around {equilibrium} $U_1$, {resulting} in
		\be
		\label{SistW}
		\displaystyle \frac{\pa{\bf W}}{\pa t}\,=\,{\mathbb D}\Delta_{\bf x}{\bf W}+{\mathbb A}{\bf W} \; \mbox{on} \; (0,\infty)\times\Gamma_{\bx},
		\ee
		{where t}he diffusion matrix $\mathbb D$ has all entries equal to zero except for elements in positions (3, 3), (3, 4), and (4, 4), that are equal to $\Phi_0(R_1)$, $ -\xi \Phi_1(R_1) R_1$ and $\delta$ respectively.
		Solutions of system \eqref{SistW} written in Fourier series are of type
		\be \label{Wgen} {\bf W}({\bf x},t)\,=\,\sum_{k}c_ke^{\lambda_k t}\,{{{\bf \overline W}_k}({\bf x})}\,,\ee
		where {the} wavenumbers $k \in \mathbb{N}$ and the eigenfunctions $\overline{\bf W}_k({\bf x})$ {define} the solution of the time--independent problem.
		Turing instability occurs when it is possible to find some wavenumber $\overline k$ such that $Re(\lambda_{\overline k})>0$, with $Re$ indicating the real part. Substituting \eqref{Wgen} in \eqref{SistW}, {we} easily verif{y} that each $\lambda_k$ is an eigenvalue of the matrix $\mathbb A-k^2\mathbb D$. A sufficient condition to have positive eigenvalues is that $det(\mathbb A-k^2\mathbb D)>0$, being 
		\be\nn
		\det({\mathbb A} - k^2 {{\mathbb D}}) \,=\,
		\frac{\, (\theta\,\mu- \eta\,)\,(R_1\,(R_1\, \Theta + \Xi) + \Xi\,\Omega)}{\mu\, \phi\, (R_1 + \Omega)}\,h(k^2),
		\ee
		with  the function of $k$
		\be\nn
		h(k^2):\,=\, k^4 \delta\,\Phi_0(R_1) + k^2 \left(\delta\,R_1\,\beta\,\phi - \Phi_1(R_1)\,R_1 \,\xi + \Phi_0(R_1)\,\tau\right) + \tau\,R_1\,\beta\,\phi.
		\ee
		Relying on { existence conditions}
		\eqref{EqPos}, the requirement $\det({\mathbb A} - k^2 {{\mathbb D}})> 0$
		is satisfied if and only if, for some values of $k$, we have $h(k^2)<0$, and this condition is equivalent to 
		\be  \nn
		\label{TurCond}
		\delta\,R_1\,\beta\,\phi  - \Phi_1(R_1)\,R_1 \,\xi + \Phi_0(R_1)\,\tau<0,\mbox{ and } h_{min} >0,\ee
		\be \nn
		h_{min} \,=\, \left(\delta\,R_1\,\beta\,\phi  - \Phi_1(R_1)\,R_1 \,\xi + \Phi_0(R_1)\,\tau\right)^2-4\,\delta\,R_1\,\beta\,\phi \,\Phi_0(R_1)\tau.
		\ee
		The two conditions above {may be} summarized {as}
		\be 
		\xi> \frac{2\sqrt{\delta\, \Phi_0(R_1)\, \tau\,R_1\,\beta\,\phi}+\delta\,R_1\,\beta\,\phi +\Phi_0(R_1)\,\tau}{\,\Phi_1(R_1)\, R_1}.
		\label{eq:pf} 
		\ee
		
		\subsection{Numerical simulations}
		
		For illustrative purposes, we {perform some numerical simulations to show the pattern formations captured by our model.
			We} set the following values for the parameters 
		\be\label{Pars}
		\beta = 0.2,\,\,\,\,\zeta = 2,\,\,\,\,\mu = 2.01,\,\,\,\,\nu = 1,\,\,\,\,
		\tau = 0.5,\,\,\,\,
		\eta = 1,\,\,\,\,
		\phi = 1,\,\,\,\,
		\delta=0.1.
		\ee
		and take the function $\varphi_1$ as {$\varphi_1 {(y)}\,=\,\cos\left(\frac{\pi}{2}\,y\right)\,\mathbf1_{[0,R_m)]}$, being $\mathbf 1_U$ the indicator function on the set $U$.}
		
		With this choice, we find the values
		\be \nn
		\bar\theta\,=\,0.49,\quad \bar\theta-\beta\,\phi\,=\,0.29,\quad \theta_-\,=\,-0.53,\quad \theta_+\,=\,0.51,
		\ee
		and we report in Figure \ref{fig1} the bifurcation diagram in the space $(\theta, \xi)$. In accordance with dashed vertical lines, representing critical values for $\theta$, we observe that the stability region lies within the range of values that guarantee the admissibility of equilibrium $U_1$ for the model. Moreover, values of $\theta$ and $\xi$  in the light-blue region satisfy the condition \eqref{eq:pf}  for the Turing instability.
		
		\begin{figure}[ht!]
			\centering
			\includegraphics[scale=0.4]{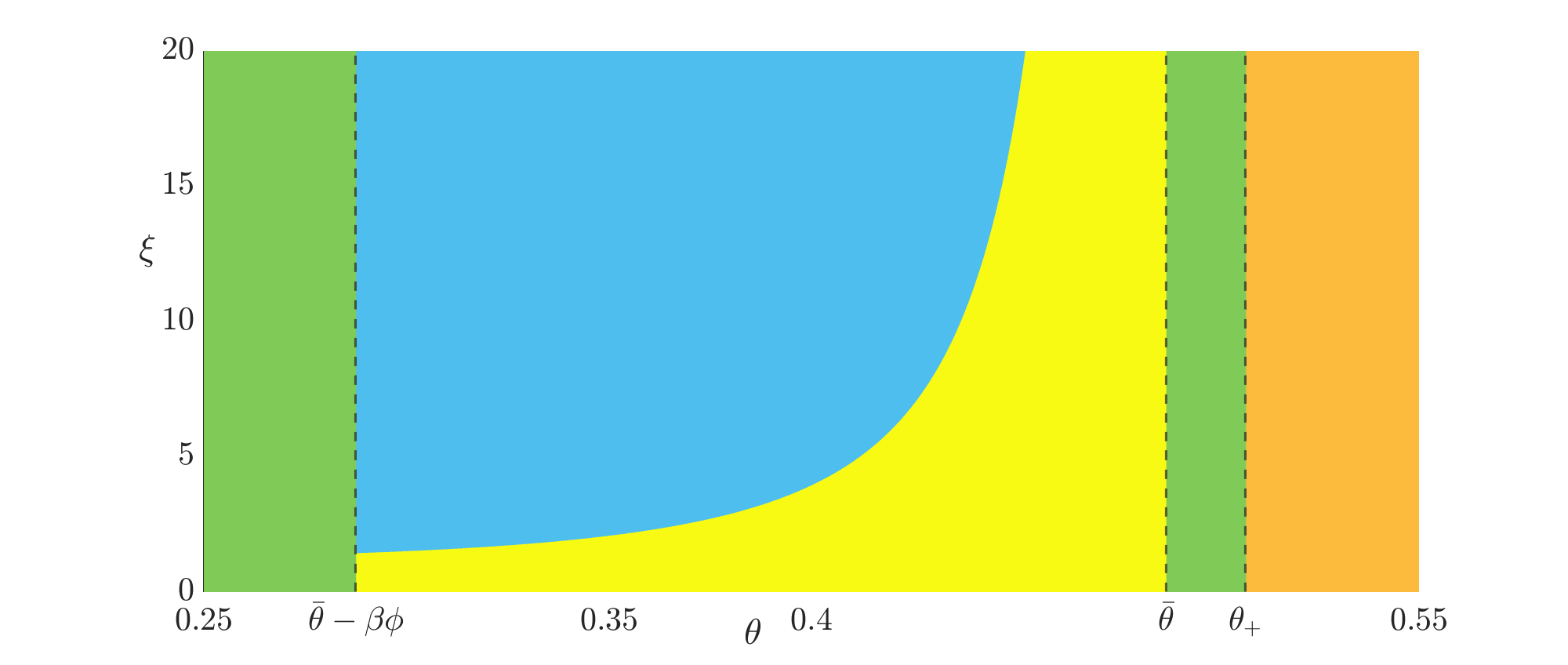}
			\caption{Bifurcation diagram in the space of parameters for system \eqref{eq:rdsA}, taking parameters as in \eqref{Pars}. Dashed lines are plotted in correspondence of values of $\overline\theta$ given in \eqref{EqPos}, $\overline\theta-\beta\,\phi$ and $\theta_+$ given in \eqref{ThetaPM}. Values in the light-blue region are those satisfying condition \eqref{eq:pf}  for the Turing instability.}
			\label{fig1}
		\end{figure}
		
		{W}e perform numerical simulations to reproduce the consumption of myelin in the white matter caused by Multiple Sclerosis. In particular, our aim is to {choose a proper set} of parameters in such a way that the formation of patterns may be compared {to} different phases of the disease, {namely in the relapsing-remitting stage (RRMS), 
			secondary progressive phase (SPMS) and primary progressive disease (PPMS).}
		Thus, we consider a one-dimensional domain of size $L\,=\,7\pi$ and we take as initial condition a random perturbation of the equilibrium state  $U_1$ for $A,\,S,\,R,\,C,$ while we initially set $E\,=\,0$. Then, we simulate the behavior of the system  \eqref{eq:rdsA} along with {the} zero-flux boundary conditions \eqref{ZeroFlux}.

		We take parameters as in \eqref{Pars} and choose the following reference parameters, {reg}arding consumption and restoration of the myelin
		\be \label{Pars2}
		{ \Theta\,=\,30,}\quad \Xi\,=\,0.02,\quad \Omega\,=\,0.001.
		\ee
		We fix the value of $\theta=0.42$, then we explore { three different cases} varying the value of $\xi$. 
		
		{ \subsubsection*{Case 1: { $\xi=6$}} We {start} by taking { low } value. The results are shown in Figure \ref{figu2}. {In} particular, we report the density of self-reactive leukocytes (Panel (a)), cytokines (Panel (b)), and myelin consumed (Panel (c)) in a time interval $0<t\leq 500$. In this case is possible to see how the model shows the following behavior:
			\begin{itemize} 
				\item Concentration of self-reactive leukocytes is driven by the higher amount of cytokines in certain areas of the domain, in which the myelin consumption is higher.
				\item Areas in which there is major demyelination (yellow spots on Panel(c)) undergo a process of remyelination (blue areas in Panel (c)), reproducing shadow plaques.
				\item This behavior {replicates over time,} describing both the acute inflammation in active plaques and the subsequent restoration performed by oligodendrocytes.
			\end{itemize}
			
			We notice that this case suitably reproduces the relapsing-remitting phase.  }

		\begin{figure}[ht!]
			\centering
			\includegraphics[scale=0.37]{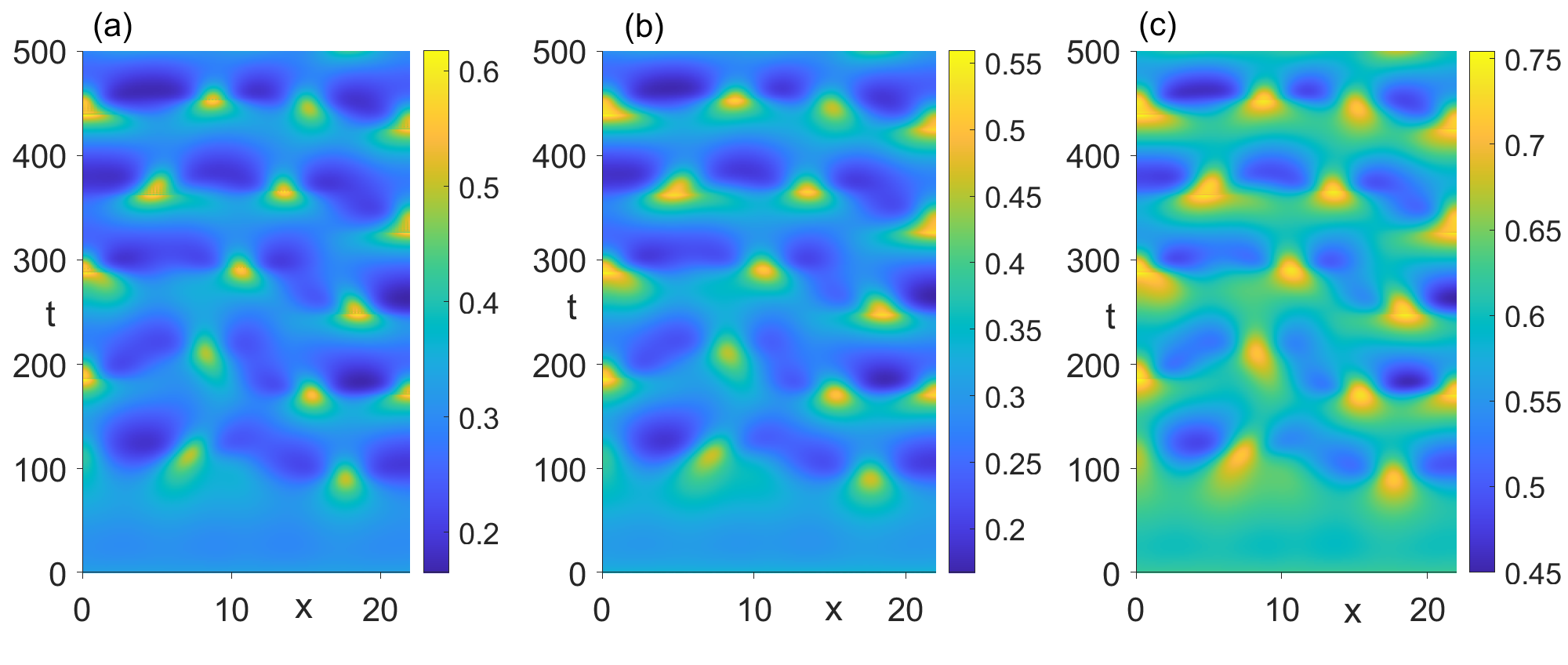}
			\caption{Behavior in space and time {typical of phase RRMS. D}ensities for self-reactive leukocytes (Panel (a)), Cytokines (Panel (b)), and destroyed myelin (Panel(c)), taking parameters as in \eqref{Pars} and \eqref{Pars2}, for values of $\theta=0.42$ and $\xi=6$.}
			\label{figu2}
		\end{figure}

		{ \subsubsection*{Case 2:  { $\xi=9$}} We increase the value of the chemotaxis parameter. In this case, as shown in Figure \ref{figu3}, in which only the behavior of destroyed myelin is reported, the choice of parameters provides the following dynamics:
			\begin{itemize} 
				\item For initial time  $0<t\leq 500$ (Panel (a)) a behavior analogous to the one obtained in the previous case can be appreciated.
				\item For a longer time $500<t\leq 2000$ (Panel (b)) a few areas of damaged myelin remain stable in time.
			\end{itemize}
			This second case accurately depicts the transition from the relapsing-remitting phase to the secondary progressive phase.}
		\begin{figure}[ht!]
			\centering
			\includegraphics[scale=0.35]{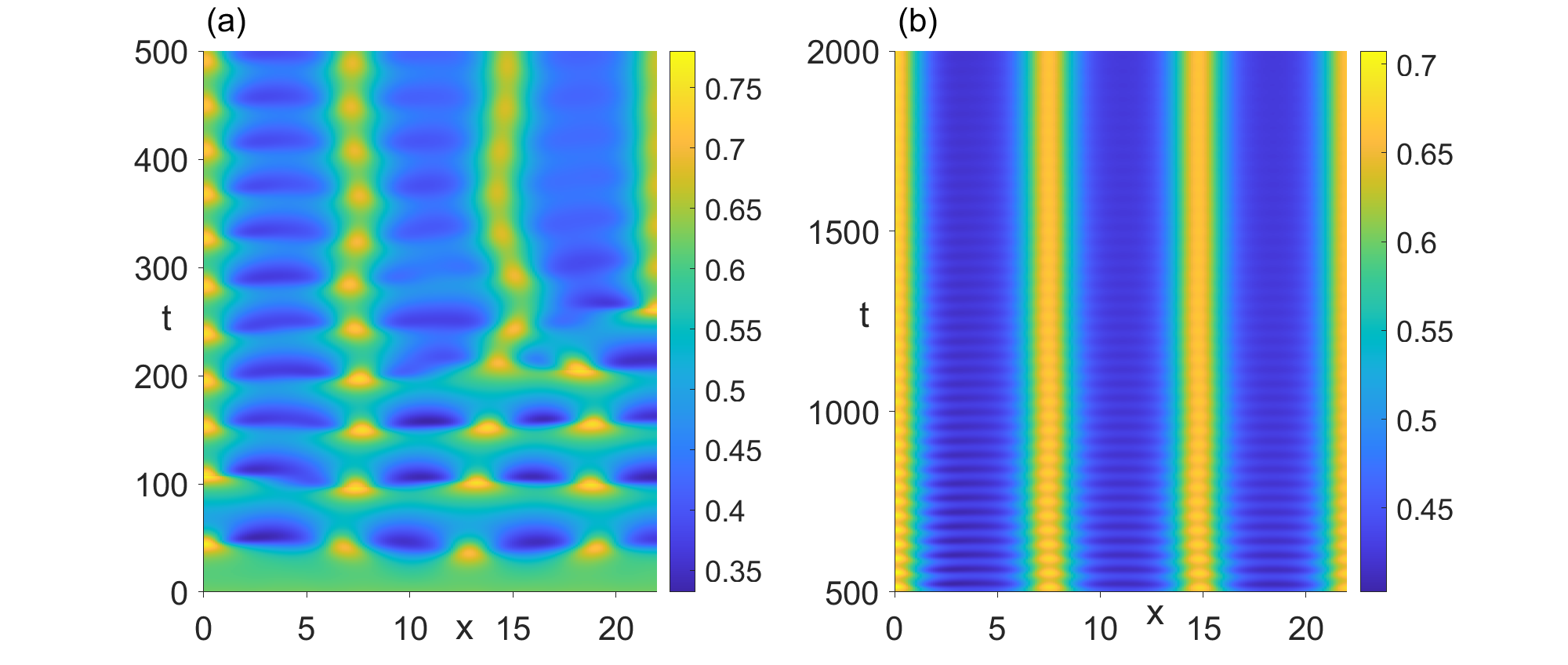}
			\caption{Behavior in space and time {typical of phase RRMS (Panel (a)) and {phase SPMS} (Panel (b)). D}densitiy for destroyed myelin in two different time intervals, taking parameters as in \eqref{Pars} and \eqref{Pars2}, for values of $\theta=0.42$ and $\xi=9$.}
			\label{figu3}
		\end{figure}

		{ \subsubsection*{Case 2: { $\xi=16.5$}} Finally, we take a higher value $\xi=16.5$. As reported in Figure \ref{figu4} for the density of self-reactive leukocytes (Panel (a)), cytokines (Panel (b)), and myelin consumed (Panel (c)), there are no initial relapsing-remitting dynamics, but areas subject to demyelination stay stable since primary formation,
			without undergoing remyelination. This behavior can be associated with the early formation of permanent lesions typical of the primary progressive phase.}
		\begin{figure}[ht!]
			\centering
			\includegraphics[scale=0.37]{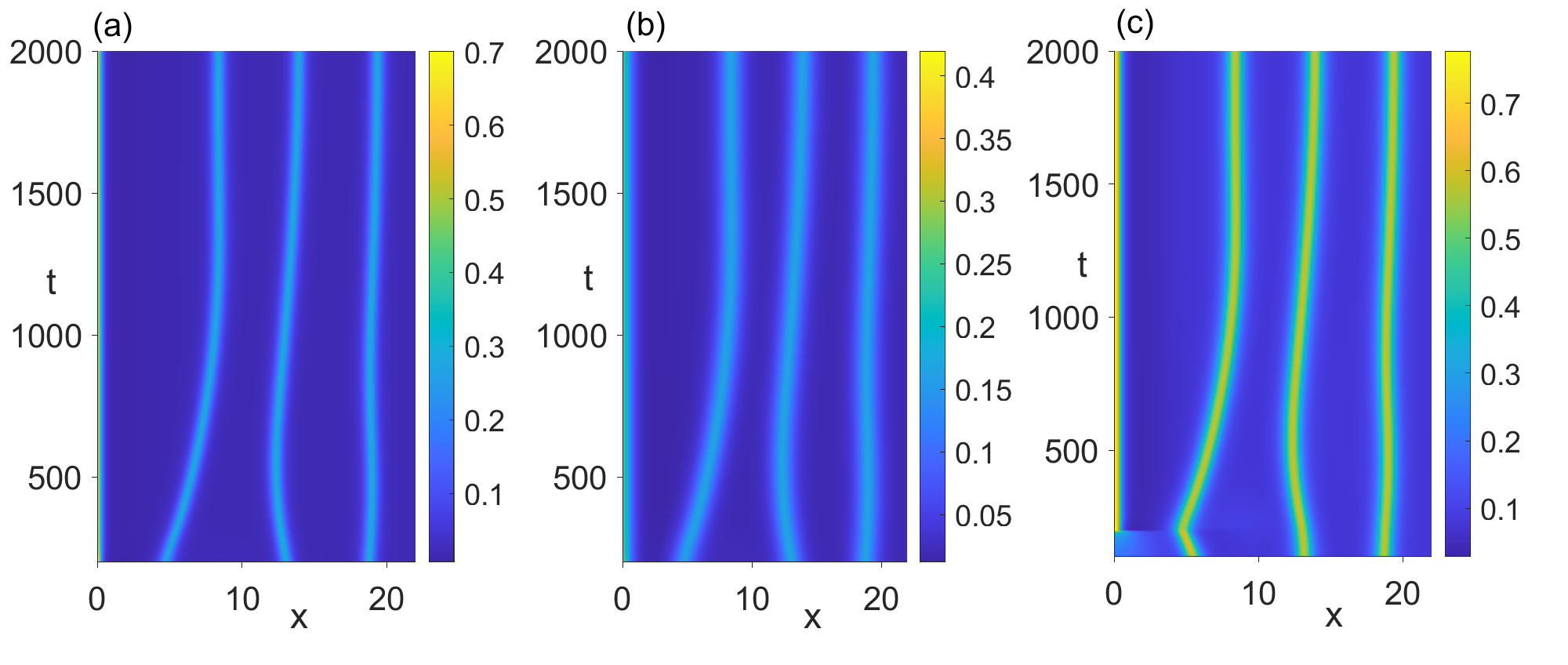}
			\caption{Behavior in space and time {typical of phase PPMS. D}ensities for self-reactive leukocytes (Panel (a)), Cytokines (Panel (b)), and destroyed myelin (Panel(c)), taking parameters as in \eqref{Pars} and \eqref{Pars2}, for values of $\theta=0.42$ and $\xi=16.5$.}
			\label{figu4}
		\end{figure}
		
		{From {the} results obtained in this section, we can appreciate the role of the chemotaxis parameter $\xi$ in
			{dynamics determined by our model}. 
			We underline {that $\xi$ is related to the microscopic parameter $\gamma$ by \eqref{DRChi} and \eqref{CoefAdim},}
			that gives the rate of velocity jumps for self-reactive leukocytes induced by pro-inflammatory cytokines gradient. 
			Specifically, as the intensity of chemotactic attraction increases, the diffusion of self-reactive leukocytes decreases, impeding remyelination and resulting in the formation of permanent lesions.}

		\section{Conclusions}
		\label{Sec6}
		
		We have presented a detailed modeling appropriate to the study of human cell dynamics, occurring in the pathological conditions of Multiple Sclerosis. The main idea behind this work was to model the destruction of myelin in the central nervous system caused by self-reactive leukocytes. For this purpose, we have included space dependence into the model presented in  \cite{della2022mathematical} (where the dynamics between cells involved in autoimmune disease were studied only with respect to time). More precisely, we considered various spatial effects like diffusion and chemotactic motion of the leukocytes driven by cytokines.
		
		We have introduced the  equations for the distribution functions of the cell populations. These include the conservative and non-conservative interaction terms, analogous to the ones in  \cite{ramos2019kinetic}. We also considered constant input of self-antigen presenting cells and natural decays of cells proposed in  \cite{della2022mathematical}. The other additional terms account for changes in the velocity of the cells and cytokines, by means of turning operators. We used them to describe the migration of the self-reactive leukocytes, regulated by the gradient of concentration of cytokines. This bias must, in turn, be distinguished from the random movement of particles,  responsible for the diffusion of both cells and cytokines.
		
		To derive the macroscopic equations, we have considered the different time scales for which each effect is relevant. The interactions that change the velocity are taken as the dominant effects at shorter time scales. The conservative and non-conservative terms, instead,  are significant only at longer time scales. Successively, we have performed a diffusive limit, obtaining a macroscopic system for the population densities. In this system, diffusion and chemotaxis play a part in the evolution of the self-reactive leukocytes and cytokines. Lastly, we have derived the equation for the destroyed myelin, that will naturally depend on the concentration of self-reactive leukocytes, while including possible remyelination processes. In this case, the division of demyelination and remyelination in stages, occurring at different time scales, leads to a nonlinear term in the equation.
		
		Finally, we have performed a Turing instability analysis of the whole system. We first outlined the stability of the system without diffusion and then we were able to find conditions on parameters for which the diffusive terms create instabilities. Instability leads to pattern formation, especially for the destruction of myelin, that is expected from Multiple Sclerosis. This result was then observed through numerical simulations. In particular, we have pointed out the relevance of the chemotactic parameter in describing the diverse phases of the diseases, as they are reported in medical literature.
		
		The present work aimed to enrich the literature on kinetic { {theory modeling}}  for autoimmune diseases, to provide some advances in the study {of Multiple Sclerosis, and, perhaps, to contribute with some outputs that can be useful in the treatment of the disease.} A more complex, two-dimensional chemotaxis model could, in principle, be used to compare numerical results with the actual lesion formation as measured in patients with Multiple Sclerosis, allowing, then, to fit the macroscopic parameters accordingly. Moreover, also the {{microscopic}} description may be refined, including the effects of treatment and drug delivery. This would give scope for future research that, if successful, may have a significant impact on patient treatment due to the consistency of a mathematical model in predicting the development of the disease over time.

		\appendix

		\section{Remarks on positivity and boundedness}
		\label{App}
		
		We start by considering the system \eqref{eq:macA} without chemotaxis and diffusion, and we rewrite the system as an equation
		\begin{equation}
			\frac{\partial \mathbf{U}(t, \bx)}{\partial t} \,=\, F\left[\mathbf{U}(t, \bx)\right] \label{eq:v},
		\end{equation}
		recalling that $\mathbf{U}(t, \bx) \,=\, (A, S, R, C, E)$ and where $F$ is a function of $\mathbf{U}$ that represents the left-hand side of equations \eqref{eq:macA}. Note that $F$ only depends on the functions that comprise the vector $\mathbf{U}$. 
		
		If we consider initial homogeneous conditions $\mathbf{U} (0, \bx ) \,=\, \mathbf{U}_0 $, it is straightforward to show that the system without diffusion and chemotaxis does not depend on the position for all time $t$. { It is easy to show that the equation \eqref{eq:v} admits solution, at least locally in time, that can be written as}
		\begin{equation}\nn
			\mathbf{U}(t, \bx) \,=\, \mathbf{U}(t_0, \bx) + \int_{t_0}^t  F\left[\mathbf{U}(t, \bx)\right]dt \ .
		\end{equation}
		We can rewrite this equation as an iterative process. Considering a time step $\Delta t$, the solution for time step $i+1$ is given by
		\begin{equation}\nn
			\mathbf{U}(t_{i+1},\bx) \approx \mathbf{U}(t_i,\bx) +  F\left[\mathbf{U}(t_i, \bx)\right] \Delta t \ .
		\end{equation}
		For $i\,=\,0$ and $t_0\,=\,0$ it results in
		\begin{equation}\nn
			\mathbf{U}(t_{1},\bx) \approx \mathbf{U}(0,\bx) +  F\left[\mathbf{U}(0, \bx)\right] \Delta t \,=\, \mathbf{U}_0 + F(\mathbf{U}_0) \ .
		\end{equation}
		We can do this iteratively for any $n$ steps and find
		\begin{equation}\nn\label{vtnh}
			\mathbf{U}(t_n,\bx) \approx \mathbf{U}_0 + F\left(\mathbf{U}_0\right)\Delta t + F\left(\mathbf{U}_0 + F\left(\mathbf{U}_0\right)\Delta t\right)\Delta t + ...
		\end{equation} 
		Note that the final result is independent of $\bx$, meaning that $\mathbf{U}$ stays homogeneous throughout the whole evolution of the system. This makes sense as the proliferation, destruction, and death of cells will be uniform throughout space as long as the cell density is homogeneous. Therefore $\mathbf{U}(t,\mathbf{x}) \,=\, \mathbf{U}(t)$ and the vector of total cell populations $\mathbf{n}(t)\,=\,(n_A(t), n_S(t), n_R(t), n_C(t), n_E(t))$, where $n_i(t)$ is the number of cells in the corresponding population $i$, can be written as:
		\begin{equation}\nn
			\mathbf{n}(t) \,=\, \int \mathbf{U} (t,\bx) d\bx \,=\, \mathbf{U}(t) \int d\bx \,=\, \mathbf{U}(t) |V| ,
		\end{equation}
		where $|V|$ is the total volume of the system.
		
		Then, {the dynamics of the number of cells can be specified if} we integrate equations \eqref{eq:macA} without chemotaxis and diffusion to obtain the following system
		\be
		\begin{cases}
			\dot{n}_A(t) \,=\, \alpha^* + p^*_{12} n_A(t) n_R(t) - d^*_{13} n_A(t) n_S(t) - d_1 n_A(t), \label{eq:MA}\\[3mm]
			\dot{n}_S(t) \,=\, p^*_{31} n_S(t) n_A(t) - d_3 n_S(t), \\[3mm]
			\dot{n}_R(t) \,=\, p^*_{21} n_R(t) n_A(t) - d^*_{23} n_R(t) n_S(t) -d_2 n_R(t),\\[3mm]
			\dot{n}_C(t) \,=\, p^*_{C2} n_A(t) n_R(t) - d_C n_C(t),\\[3mm]
			\dot{n}_E(t) \,=\,  (\hat{n}_E - n_E(t) )\dfrac{b^*_{62}\,b^*_{52}n_R(t)}{r_5+b^*_{52}n_R(t)}n_R(t) - r_6 n_E(t),
		\end{cases}
		\ee
		where the asterisk-labeled coefficients corresponds to a rescaled quantity $p^* \,=\, p/|V|$ and $\hat{n}_E\,=\,\bar{n}_E |V|$ due to the size of the system. Boundedness and positivity of the system \eqref{eq:MA} was shown in the paper\cite{della2022mathematical} for $p_{21}^* < p_{31}^*$ and we comment here on the more general system that includes last two equations of \eqref{eq:MA}. The solution of the equation for $C$ is
		\be \nn
		n_C(t) \,=\, n_C(0) + e^{-d_C t}\int_0^t e^{d_C s} n_A(t)n_R(t) ds \ .
		\ee
		Both boundedness and positivity of $n_C(t)$ for all $t$ come straightforwardly from the boundedness and positivity of both $n_A(t)$ and $n_R(t)$, for positive initial value $n_C(0)$.
		
		As for $n_E(t)$ we obtain
		\be \nn
		n_E(t) \,=\, e^{-I(t)}\left[n_E(0) + \hat n_E\int_0^t e^{I(u)}\frac{b^*_{62}\,b^*_{52}n_R(u)}{r_5+b^*_{52}n_R(u)}n_R(u) du\right] ,
		\ee
		where
		\be \nn
		I(t) \,=\, \int_0^t\left(\frac{b^*_{62}\,b^*_{52}n_R(t)}{r_5+b^*_{52}n_R(t)}n_R(t)+r_6\right)ds \ .
		\ee
		Naturally, as long as $n_R(t)$ is positive and bounded, the integral $I(t)$ will also be positive and bounded (assuming positive coefficients) and, as a consequence, so will $n_E(t)$ for positive initial value $n_E(0)$.
		For the inhomogeneous case, a similar reasoning as the homogeneous case for every point in space can be used. Due to the lack of spatial derivatives in the system, the value of $\mathbf{U}(t,\mathbf{x_0})$ for any fixed point $x_0$ will only depend on the initial value of that point $x_0$, as there is no flux of cells between points. Rewriting \eqref{vtnh} for the inhomogeneous case for any fixed point $x_0$, and denoting $\mathbf{U}(0,\bx_0) \,=\, \mathbf{U}_0(\bx_0)$, we get:
		\be \nn
		\mathbf{U}(t_n,\bx_0) \approx \mathbf{U}_0(\bx_0) + F\big[\mathbf{U}_0(\bx_0)\big]\Delta t + F\bigg[\mathbf{U}_0(\bx_0) + F\big[\mathbf{U}_0(\bx_0)\big]\Delta t\bigg]\Delta t + ...\, .
		\ee
		As we can see, $\mathbf{U}(t_n,\bx_0)$ will only depend on $\mathbf{U}_0(\bx_0)$ and the time intervals $\Delta t$. So we can construct a time-independent system for every point $\bx_0$ just like the system \eqref{eq:MA} with different coefficients. The positivity and boundedness follow in the same way as the homogeneous case.
		
		\section*{Acknowledgment}
		The research of both authors has been carried out in the frame of activities sponsored by the Cost Action CA18232. Author RT also thanks the GNFM (National Group of Mathematical-Physics) of INdAM (National Institute of Advanced Mathematics).  The work of both authors was supported by the Portuguese FCT/MCTES  (Fundação para a Ciência e a Tecnologia) Projects  UIDB/00013/2020,  UIDP/00013/2020 and PTDC/03091/2022 (“Mathematical Modelling of Multi-scale Control Systems: applications to human diseases (CoSysM3)”).

\bibliographystyle{abbrv}
\bibliography{biblioAbb}

\end{document}